\documentclass[11pt]{article}
\usepackage[american]{babel}
\usepackage{amsmath}
\usepackage{latexsym}
\usepackage{amssymb}
\usepackage{theorem}
\usepackage{diagrams}
\usepackage{pstricks,pst-3d}
\newgray{gray9}{.9}
\newgray{gray8}{.8}
\newgray{gray7}{.7}
\newgray{gray6}{.6}
\newgray{gray5}{.5}
\newgray{gray4}{.4}
\newgray{gray3}{.3}
\newgray{gray2}{.2}
\newgray{gray1}{.1}
\usepackage{calligra}
\DeclareMathAlphabet{\mathcalligra}{T1}{calligra}{m}{n}
\DeclareFontShape{T1}{calligra}{m}{n}{<->s*[1.2]callig15}{}
\theoremstyle{change}
\newtheorem{Thm}{Theorem}[section]
\newtheorem{Cor}[Thm]{Corollary}
\newtheorem{Prop}[Thm]{Proposition}
\newtheorem{Lem}[Thm]{Lemma}
{\theorembodyfont{\rmfamily}
\newtheorem{Num}[Thm]{}

\newtheorem{Def}[Thm]{Definition}}
\renewcommand{\phi}{\varphi}
\newcommand{\bra}[1]{\langle#1\rangle}
\newcommand{\proof}{\par\medskip\rm\emph{Proof. }}
\newcommand{\skop}{\par\medskip\rm\emph{Sketch of proof. }}
\newcommand{\qed}{\ \hglue 0pt plus 1filll $\Box$}
\newcommand{\mapstoo}{\longmapsto}
\newcommand{\RR}{\mathbb{R}}
\newcommand{\ZZ}{\mathbb{Z}}
\newcommand{\NN}{\mathbb{N}}
\renewcommand{\SS}{\mathbb{S}}

\newcommand{\CC}{\mathbb{C}}
\newcommand{\id}{\mathrm{id}}
\newcommand{\SKIP}[1]{}
\newcommand{\eps}{\varepsilon}
\renewcommand{\emptyset}{\varnothing}
\newcommand{\proj}{\mathrm{proj}}
\newcommand{\F}{\mathcal F}

\newcommand{\partialmax}{{\partial}_{\mathit{cpl}}}
\newcommand{\A}{\mathcal A}
\newcommand{\Amax}{\hat{\mathcal A}}
\newcommand{\of}{\mathrm{of}}
\newcommand{\ms}{\mathrm{ms}}

\newcommand{\fin}{\mathrm{fin}}
\newcommand{\infi}{\mathrm{inf}}
\newcommand{\Res}{\mathrm{Res}}
\newcommand{\Opp}{\mathrm{Opp}}
\newcommand{\EB}{\mathrm{EB}}
\newcommand{\Isom}{\mathrm{Isom}}
\renewcommand{\setminus}{-}
\begin{document}
\title{{\bf Coarse Equivalences of Euclidean Buildings}}
\author{Linus Kramer and Richard M.~Weiss%
\footnote{Kramer and Weiss were
supported by the SFB 878 and DFG Project KR 1668/7.}
\\ \\With an appendix by Jeroen Schillewaert and Koen Struyve}
\date{}
\maketitle
\begin{abstract}
We prove the following rigidity results.
Coarse equivalences between metrically complete
Euclidean buildings
preserve spherical buildings at infinity.
If all irreducible factors have dimension at least two, then
coarsely equivalent
Euclidean buildings are isometric (up to scaling factors);
if in addition
none of the irreducible factors is a Euclidean cone, then the
isometry is unique and has finite distance from the
coarse equivalence. 

The appendix shows how these results can be extended to non-complete Euclidean
buildings.
\end{abstract}

We prove coarse (i.e. quasi-isometric) rigidity results
for leafless trees (simplicial trees and $\RR$-trees with extensible geodesics) 
and, more generally, for
discrete and nondiscrete Euclidean buildings. 
For trees, a key ingredient is a certain equivariance condition. 
Our main results are as follows.

\medskip\noindent\textbf{Theorem I\ }
{\em Let $G$ be a group acting isometrically on two metrically complete leafless
trees $T_1,T_2$. Assume that there is a coarse equivalence
$f:T_1\rTo T_2$, that $T_1$ has at least $3$ ends
and that the induced map $\partial f:\partial T_1\rTo\partial T_2$
between the ends of the trees is $G$-equivariant. If the
$G$-action on $\partial T_1$ is $2$-transitive, then (after rescaling the
metric on $T_2$) there is a 
$G$-equivariant isometry $\bar f:T_1\rTo T_2$
with $\partial f=\partial \bar f$. If $T_1$ has at least two branch points,
then $\bar f$ is unique and has finite distance from $f$.}

\medskip\noindent
The precise result is \ref{MainTreeProp} below,
where we also consider products of trees and Euclidean spaces. 
Without the equivariance condition, quasi-isometric rigidity fails.
The letters \textsf{X} and \textsf{H} written infinitely large are, for example, trees
which are coarsely equivalent without
being isometric. Note, too, that any two simplicial trees of finite constant
valence are coarsely equivalent without being necessarily isometric \cite{Papa}.

Our next results are concerned with (possibly nondiscrete) Euclidean buildings.

\medskip\noindent\textbf{Theorem II\ }
{\em Let $X_1$ and $X_2$ be Euclidean
buildings whose spherical buildings at infinity $\partialmax X_1$
and $\partialmax X_2$ are thick.
Suppose that $f:X_1\times\RR^{m_1}\rTo X_2\times\RR^{m_2}$ is a coarse equivalence.
Then $m_1=m_2$ and there is a combinatorial isomorphism
$f_*:\partialmax X_1\rTo\partialmax X_2$ between the
spherical buildings at infinity which is characterized by the fact
that for an affine apartment $A\subseteq X_1$ the $f$-image of 
$A\times\RR^{m_1}$
has finite Hausdorff distance from $f_*(A)\times\RR^{m_2}$.}

\medskip\noindent
This is \ref{MainThmII} below.
We remark that the boundary map $f_*$ is constructed in a combinatorial
way from $f$. In general, a coarse equivalence between CAT$(0)$-spaces
will not induce a map between the respective Tits boundaries. 

\medskip\noindent\textbf{Theorem III\ }
{\em Let $f:X_1\times\RR^m\rTo X_2\times\RR^m$ be as in Theorem II and
assume in addition that $X_1$ has no tree factors and that $X_1$ and
$X_2$ are metrically complete. 
Then there is (possibly after rescaling the metrics
on the de Rham factors of $X_2$) an isometry
$\bar f:X_1\rTo X_2$ with $(\bar f\times\id_{\RR^m})_*=f_*$.
Put $f(x\times v)=f_1(x\times v)\times f_2(x\times v)$.
If none of the de Rham factors of $X_1$ is a Euclidean cone over
its boundary, then $\bar f$ is unique and $d(f_1(x\times y),\bar f(x))$ is
bounded as a function of $x\in X_1$.}

\medskip\noindent
For a more general statement see \ref{MainTheoremOfPaper} below.

Kleiner and Leeb proved results similar to Theorems II and III in 
\cite[1.1.3]{KL} and \cite[1.3]{Le}
under the additional hypothesis that the Euclidean buildings are
simplicial and locally finite or 
their buildings at infinity satisfy the Moufang condition.
Their work extended Mostow-Prasad rigidity \cite{Pra}. Our results offer
in particular an alternative approach to these important achievements.

\medskip\noindent
\textbf{Outline of the paper\ }
In the first section we collect some basic facts about metric spaces,
nonpositive curvature, ultraproducts and ultralimits.
The second section is concerned with
trees. We derive in particular the structure of leafless trees that admit an
isometry group that is $2$-transitive on the ends. 
Theorem I is then an application of this  structural result, combined with
Morse's Lemma. In the third section we collect the facts about spherical
buildings and projectivities that we need. The fourth section is devoted
to Euclidean buildings and their local and global structure.
In the fifth section we combine all previous results. We first prove
Theorem II, using the `higher dimensional Morse Lemma'. Then we
combine Theorem I and Theorem II in order to prove Theorem III.
The proof relies, among other things, on  Tits' rigidity result 
\cite[Thm.~2]{TiComo} of `ecological' (tree-preserving)
boundary isomorphisms of Euclidean buildings.

\section{Some metric geometry}
We recall a few notions and facts about metric spaces that will be needed later.
Let $(X,d)$ be a metric space. For $r>0$ and $x\in X$ we put 
\[
B_r(x)=\{y\in X\mid d(x,y)<r\}
\quad\text{and}\quad
\bar B_r(x)=\{y\in X\mid d(x,y)\leq r\}.
\]
For a subset $Y\subseteq X$ we put $B_r(Y)=\bigcup\{B_r(y)\mid y\in Y\}$.
We call $Y\subseteq X$ \emph{bounded}
if $Y$ is contained in some sufficiently large ball, $Y\subseteq B_r(x)$ for some $x\in X$,
and we call $Y$ \emph{cobounded} if $X= B_r(Y)$ for some $r\geq 0$.

Let $f:X\rTo Y$ be a map between metric spaces. If there is a constant $r>0$ such
that 
\[
d(f(u),f(v))\leq rd(u,v)
\]
holds for all $u,v\in X$, we call $f$ an
\emph{$r$-Lipschitz map}. A map  $f$ which preserves distances
is called an \emph{isometric embedding}; if $f$ is in addition onto, it is called an \emph{isometry}.

The group of all isometries of $X$ onto itself is the \emph{isometry group} $\Isom(X)$.
An \emph{isometric action} of a group $G$ on $X$ is a homomorphism $G\rTo \Isom(X)$.
We call such an action (or group) \emph{bounded} if some (or, equivalently, every)
$G$-orbit is bounded in $X$. The action is called \emph{cobounded} if there is a
cobounded orbit.

\begin{Num}\textbf{\boldmath CAT$(\kappa)$ spaces\ }
\label{CAT(k)Definition}
A \emph{geodesic} in a metric space $X$ is
an isometric embedding $\gamma:J\rTo X$, where $J\subseteq \RR$ is a
closed interval.
The image $\gamma(J)$ is then called a \emph{geodesic segment}.
If the domain of $\gamma$ is $J=\RR$, then $\gamma(\RR)$ is also called
a \emph{geodesic line}, and if $J=[0,\infty)$, then $\gamma([0,\infty))$
is called a \emph{geodesic ray.} 
If any two points of $X$ are contained in some
geodesic segment, $X$ is called a \emph{geodesic space}. If every geodesic
$\gamma:J\rTo X$ admits a geodesic extension $\bar\gamma:\RR\rTo X$, we say that
$X$ has \emph{extensible geodesics}.

A geodesic space $X$ is called a \emph{CAT(0) space} if no geodesic triangle
in $X$ is thicker than its comparison triangle in Euclidean space $\RR^2$,
see \cite[II.1]{BH}. More generally, a geodesic space is called
a CAT$(\kappa)$ space, for $\kappa\in\RR$, if no geodesic triangle in
$X$ is thicker than its comparison triangle in the complete simply connected
Riemannian $2$-manifold $M_\kappa$ of constant sectional curvature $\kappa$.
For $\kappa>0$, the space $M_\kappa$ is a sphere and the condition
is only required for geodesic triangles of perimeter less than 
$2\frac\pi{\sqrt\kappa}$ (and the existence of geodesics is only required
between points at distance less than $\frac\pi{\sqrt\kappa}$).
If $(X,d)$ is CAT$(\kappa)$, then the same space is also CAT$(\kappa')$
for all $\kappa'\geq\kappa$ \cite[II.1.12]{BH}. If $r>0$ and if $(X,d)$ is CAT$(\kappa)$,
then the rescaled space $(X,rd)$ is CAT$(\kappa/{\sqrt r})$.

If $K$ is a metrically complete convex subset in a CAT$(0)$ space $X$,
then there is a $1$-Lipschitz retraction 
\[
\pi_K:X\rTo K
\]
which maps $x\in X$ to the unique closest point in $K$ \cite[II.2.4]{BH}.
\end{Num}
The metric completion of a CAT$(\kappa)$ space is again a CAT$(\kappa)$ space.
One important fact about complete CAT$(0)$ spaces is the Bruhat-Tits
Fixed Point Theorem: every bounded isometry group in such a space has
a fixed point \cite[II.2.8]{BH}.
\begin{Num}\textbf{Controlled maps\ }
We now recall some notions from coarse geometry \cite{Roe}.
A map $f:X\rTo Y$ between metric spaces is called
\emph{controlled} if there is a monotonic real function
$\rho:\RR_{\geq 0}\rTo\RR_{\geq 0}$ such that
\[
d_Y(f(x),f(y))\leq\rho(d_X(x,y))
\]
holds for all $x,y\in X$. 
If in addition the preimage of every
bounded set is bounded, then $f$ is called a \emph{coarse map}.
Neither $f$ nor $\rho$ is required to be continuous.
Note that the image of a bounded set under a controlled map is
bounded.
Two maps $g,f:X\pile{\rTo \\ \rTo} Y$ between metric spaces have \emph{finite distance}
if the set $\{d_Y(f(x),g(x))\mid x\in X\}$ is bounded. This is an
equivalence relation which
leads to the \emph{coarse metric category} whose
objects are metric
spaces and whose morphisms are equivalence classes of coarse maps.
A \emph{coarse equivalence} is an isomorphism in this category.
\end{Num}
\begin{Lem}
If $f:X\rTo Y$ and $g:Y\rTo X$ are controlled and if
$g\circ f$ has bounded distance from the identity map on $X$,
then $f$ is coarse.
In particular, $f$ is a coarse equivalence if $f\circ g$ also
has bounded distance from the identity.

\proof
Suppose that $f(Z)$ is bounded. Since $g$ is controlled,
$g(f(Z))$ is also bounded, $g(f(Z))\subseteq B_r(x)$ for some
$x\in X$.
Now $Z$ is contained in some $s$-neighborhood of $g(f(Z))$,
so $Z\subseteq B_{r+s}(x)$ is bounded.
\qed
\end{Lem}
\begin{Num}\textbf{Quasi-isometries\ }
\label{QIDef}
A controlled map $f$ with control function
$\rho(t)=t$ is the same as a $1$-Lipschitz map.
If $f$ admits a control function
\[
\rho(t)=ct+d,\qquad\text{  with }c\geq 1\text{ and }d\geq 0,
\]
then $f$ is called \emph{large-scale Lipschitz}.
A coarse equivalence is called a \emph{quasi-isometry} if both the map and
its coarse inverse are large-scale Lipschitz. 
More generally, 
we call $f$ a \emph{quasi-isometric embedding} if $f$ is a quasi-isometry
between $X$ and $f(X)\subseteq Y$.
A \emph{rough isometry} is a coarse equivalence where the control functions
in both directions are of the form $\rho(t)=t+d$ (such maps are sometimes
called \emph{coarse isometries}). A map which has finite
distance from an isometry is an example of a rough isometry.
\end{Num}
In the class of geodesic spaces, controlled maps are essentially the same as
large-scale Lipschitz maps.
\begin{Lem}
Let $f:X\rTo Y$ be controlled. If $X$ is geodesic, then
$f$ is large-scale Lipschitz. In particular, every
coarse equivalence between geodesic spaces is automatically a quasi-isometry.

\proof
For $x,y\in X$ let $d(x,y)=m+s$, with $m\in\NN$ and $0\leq s<1$.
Let $\gamma$ be a geodesic from $x$ to $y$. Then
$d(f(\gamma(j)),f(\gamma(j+1)))\leq\rho(1)$ for $j=0,\ldots,m-1$, so
$d(f(\gamma(0)),f(\gamma(m+s)))\leq m\rho(1)+\rho(s)\leq(m+s)\rho(1)+\rho(1)$.
\qed
\end{Lem}
\begin{Num}\textbf{Hausdorff distance\ }
\label{FinHd}
Two (nonempty) subsets $U,V$ of a metric space $X$ have
\emph{Hausdorff distance at most $r$} if
\[
U\subseteq B_r(V)\quad\text{ and }\quad V\subseteq B_r(U).
\]
In this case we write $Hd(U,V)<r$. We define for $U,V,W\subseteq X$
the Hausdorff distance \cite[VIII \S6]{Haus} as
\[
Hd(U,V)=\inf\{r>0\mid Hd(U,V)<r\}
\quad\text{ and }\quad
Hd_W(U,V)=Hd(U\cap W,V\cap W).
\]
For example, a nonempty
subset is bounded if and only if it has finite Hausdorff distance
from some point. 
More generally, we say that $V$ \emph{dominates} $U$ if
$U\subseteq B_r(V)$ for some $r>0$, and we write
then 
\[
U\subseteq_{Hd} V.
\]
This defines a preorder on the
subsets of $X$. 
If $f:X\rTo Y$ is controlled with control function $\rho$,
and if $U\subseteq B_r(V)$, then
$f(U)\subseteq B_{\rho(r)}(V)$. So if $f$ is a coarse
equivalence, then $U$ and $V$ have finite Hausdorff distance
if and only if $f(U)$ and $f(V)$ have finite Hausdorff
distance, and $U$ dominates $V$ if and only if $f(U)$
dominates $f(V)$.
\end{Num}
Now we turn to ultralimits and asymptotic cones of metric spaces.
These constructions generalize pointed Gromov-Hausdorff convergence
of metric spaces. The advantage of ultralimits is that they
exist even if Gromov-Hausdorff convergence fails, and that these
spaces always have good metric and functorial properties.
In contrast to \cite[Sec.~3.]{KaLe}, \cite[2.4]{KL}, \cite[I.5]{BH}, \cite[7.5]{Roe},
\cite[3.29]{Gr}, we follow the original approach 
of van den Dries and Wilkie \cite{VW}, which is
based on elementary model theory and ultraproducts. The
advantage of this is viewpoint is that it keeps track of
the global structure of spaces, and that it allows us in a natural
way to include some further structure, such as metric balls,
apartments, geodesics, or Hausdorff distance. 
We remark that similar techniques have been used in Banach
space theory for the last 50 years.

In order to make the paper self-contained, we briefly review
some elementary facts about languages, structures, and ultraproducts. 
They can be found in any `old-fashioned' textbook on mathematical logic,
such as \cite{BeSl} or \cite{CK}.
\begin{Num}\textbf{Languages and structures\ }
\label{LanguagesAndStructures}
Suppose that $\mathcal L$ is a first-order language.
Such a language is a set consisting of symbols for finitary functions, relations and constants.
The language is allowed to be infinite.
A specific example which will be important here
is the language $\mathcal L_\of=\{+,\cdot,\leq,0,1\}$
of ordered fields. An \emph{$\mathcal L$-structure} is a tuple
$\mathfrak A=(A,F,R,C)$ consisting of a nonempty set $A$ (the \emph{universe} of
the structure), a set $F$ of finitary functions, a set $R$ of finitary relations,
and a set $C\subseteq A$ of constants. 
In addition, we have an \emph{interpretation} which assigns to every
function/relation/constant symbol in $\mathcal L$ a function/relation/constant 
in $\mathfrak A$.
For example, the reals form an $\mathcal L_\of$-structure $\mathfrak R$ if
we interpret the symbols $+,\cdot,0,1,\leq$ in their usual way.
Abusing notation, we denote symbols and their interpretations by the same
letters---this should cause no confusion, as the interpretation will
always be the natural one.

Metric spaces fit well into this concept. The relevant language 
is $\mathcal L_\ms=\mathcal L_\of\cup\{d,X,R\}$.
If $(X,d)$ is a metric space, then we consider the 
$\mathcal L_\ms$-structure $\mathfrak X$ with universe
$X\stackrel{.}{\cup}\RR$, functions $+,\cdot,d$,
relations $R,X,\leq$, and constants $0,1$.
The unary relation $X(p)$ says that $p$ is a point in $X$ and
the unary relation $R(t)$ says that $t$ is a real number.
The binary function $d$ gives the distance on points in $X$,
while the binary functions $+,\cdot$ give the usual arithmetic
operations in the real numbers.
\end{Num}
\begin{Num}\textbf{Formulas\ }
From a language $\mathcal L$ we build formulas by the standard rules
of first-order logic, using the logical symbols $\exists,\forall,\dot=,\neg,\to$
and an infinite set of symbols for variables. Variables in a formula which
are not bound by a quantifier are called \emph{free}. We write
$\phi(x)$ for such a formula with a free variable $x$. Abusing notation,
$x$ may also denote a finite tuple of free variables (and we also
allow tacitly that the formula $\phi(x)$ contains no free variables at all).
If we want to interpret a formula in a given $\mathcal L$-structure, we
first have to substitute the free variables by elements of the universe, that is,
we plug in specific values $a$ for the free variables $x$.
For this substitution, we write $\phi(a)$.
A formula without free variables is called a \emph{sentence}.
If a sentence $\phi$ holds in the
$\mathcal L$-structure $\mathfrak A$, one writes
$\mathfrak A\models\phi$. For example, we have in the 
metric space $\mathfrak X$ that
\[
\mathfrak X\models\forall u,v,w\,\bigl(X(u)\wedge X(v)\wedge X(w)\to
d(x,z)\leq d(x,y)+d(y,z))\bigr)
\]
that is, the triangle inequality holds for all triples of points in $X$.
The formula 
\[
\phi(x_1,x_2)=\bigl(X(x_1)\wedge X(x_2)\wedge (d(x_1,x_2)\leq 1)\bigr)
\]
with free variables $x_1,x_2$ holds in the metric space $X$ if we
substitute points $u,v$ for $x_1$ and $x_2$ whose distance is at most $1$.
\end{Num}
\begin{Num}\textbf{Ultraproducts\ }
Suppose that $K$ is a countably infinite set and that
$D$ is a nonprincipal ultrafilter on $K$. This means that
$D$ is a collection of subsets of $K$, containing all
cofinite subsets of $K$, which is closed under
finite intersections, closed under going up (if $a\in D$ and 
$K\supseteq b\supseteq a$,
then $b\in D$), and which contains for every subset $a\subseteq K$ either
$a$ or its complement $K\setminus a$, but not both. 
Using Zorn's Lemma, it is easy to see that nonprincipal ultrafilters exist
\cite[4.1.3]{CK}.
Suppose that $(\mathfrak A_k)_{k\in K}$ is a family of
$\mathcal L$-structures $\mathfrak A_k=(A_k,F_k,R_k,C_k)$.
We introduce an equivalence relation on
$K$-sequences $(a_k),(b_k)\in\prod_kA_k$ by putting
\[
(a_k)\sim_D (b_k)\quad\text{ if }\quad\{j\in K\mid a_j=b_j\}\in D.
\]
The resulting set of equivalence classes is the \emph{ultraproduct}
$A_D=\prod_kA_k/D$. On the functions, relations and constants in 
$\prod_kF_k$, $\prod_kR_k$ and $\prod_k C_k$ we introduce the
same equivalence relation. The resulting $\mathcal L$-structure
\[
\mathfrak A_D=\bigl(A_D,F_D,R_D,C_D\bigr)
\]
is the \emph{ultraproduct} of the family $(\mathfrak A_k)_{k\in K}$.
If all structures in the family are equal, one calls $\mathfrak A_D$ also
an \emph{ultrapower}.
\end{Num}
Ultraproducts are useful because of the following two facts.
Firstly, \L os' Theorem says that $\mathfrak A_D$ has the
same first-order properties as `the majority of the $\mathfrak A_k$'.
The second fact is that the ultraproduct is `huge' and contains `strange'
elements: it is $\omega_1$-saturated  (this is essentially the same
as the `overspill' mentioned in \cite[Sec.~3]{VW}).
For the next two theorems we assume that $(\mathfrak A_k)_{k\in K}$ is
a countable family of $\mathcal L$-structures and that
$D$ is a nonprincipal ultrafilter on $K$.
\begin{Thm}[\L os' Theorem]\label{LosThm}
Let $\phi(x)$ be an $\mathcal L$-formula (possibly with a free variable $x$).
Let $a=(a_l)\in\prod_kA_k$. Then 
\[
\mathfrak A_D\models\phi(a)\quad\text{ if and only if }\quad
\{j\in K\mid \mathfrak A_j\models \phi(a_j)\}\in D.
\]
(If there is no free variable in $\phi$ this means that $\phi$ holds in
$\mathfrak A_D$ if and only if the set of indices $k$ for which $\phi$ holds
in $\mathfrak A_k$ is in $D$.)
\proof
The proof is an easy induction on the complexity of formulas
in \cite[4.1.9]{CK} \cite[\S5 2.1]{BeSl} \cite[9.5.1]{Ho}.
\qed
\end{Thm}
The following fact is an elementary consequence.
\begin{Lem}
\label{UltrapowerOfFiniteStructure}
If $\mathfrak A$ is a finite structure, then the ultrapower $\mathfrak A_D$ is
isomorphic to $\mathfrak A$.
\qed
\end{Lem}
If we put $\mathfrak A_k=\mathfrak R=(\RR,+,\cdot,\leq,0,1)$, we end up with the
field of \emph{nonstandard reals} 
\[
\mathfrak R_D=({}^\ast\RR,+,\cdot,\leq,0,1).
\]
By \L os' Theorem, this is a real closed field.
If each $\mathfrak X_k$ is a metric space, then $\mathfrak X_D$
is a generalized metric space, where the distance function takes
it values in the nonstandard reals $^*\RR$. We will come back to these
generalized metric spaces below. Variants of the following result are folklore. 
For the sake of completeness, we include a proof.
\begin{Thm}[\boldmath $\omega_1$-saturation of ultraproducts]\label{overspill}
Let $\Phi$ be countable set of $\mathcal L$-formulas
in the free variable $x$.
Suppose that for every finite subset $\Psi\subseteq\Phi$
there exists an $K$-sequence $(a_k)$ such that $\mathfrak A_D\models\psi(a)$
holds for all $\psi(x)\in\Psi$. Then there exists
an $K$-sequence $(a_k)$ such that $\mathfrak A_D\models \phi(a)$ holds
for all $\phi(x)\in\Phi$.

\proof
We may assume that $K=\NN$ and that $\Phi$ is countably
infinite. Accordingly, let $\Phi=\{\phi_k(x)\mid k\in\NN\}$ be an
enumeration of $\Phi$. We have by \L os' Theorem \ref{LosThm} that
\[
 K_n=\{k\in\NN\mid k\geq n\text{ and }
 \mathfrak A_k\models\exists x\,\phi_0(x)\wedge\cdots\wedge\phi_n(x)\}\in D.
\]
Since $\bigcap_{n\geq 0}K_n=\emptyset$, there exists for each $i\in K_0$
a largest $n(i)$ with $i\in K_{n(i)}$. 
For $k\in \NN\setminus K_0$ we choose $a_k\in A_k$ arbitrarily. 
For $k\in K_0$ we choose $a_k\in A_k$ in such a way that
$\mathfrak A_k\models\phi_0(a_k)\wedge\cdots\wedge\phi_{n(k)}(a_k)$ holds.
We claim that the $K$-sequence $a=(a_k)$ has the required property.
Let $n\in\NN$.
For all $k\in K_n$ we have that $n\leq n(k)$, whence
$\mathfrak A_k\models \phi_n(a_k)$. Thus
$\mathfrak A_D\models\phi_n(a)$ holds.
\qed
\end{Thm}
As an illustration, consider the countable set of $\mathcal L_\of$-formulas
\[
\phi_n(x)=\bigl(x>0+\underbrace{1+1+\cdots+1}_{n\text{ summands}}\bigr).
\]
By \ref{overspill}, there are elements in $^*\RR$ which satisfy these formulas
simultaneously.
This means that there are elements in $^\ast\RR$ which are bigger than
every natural number. In other words, $^\ast\RR$ is a 
\emph{non-archimedean real closed field}. We call an element $r\in{}^\ast\RR$ \emph{finite}
if $|r|\leq n$ holds for some $n\in\NN$. It is readily verified that
the finite elements form a valuation ring $^\ast\RR_\fin\subseteq{}^\ast\RR$. 
Its maximal ideal is the set $^\ast\RR_\infi=\{r\in{}^\ast\RR\mid
|r|\leq \frac1{n+1}\text{ for all }n\in\NN\}$ of \emph{infinitesimal elements}.
One can show that the quotient field $^\ast\RR_\fin/^\ast\RR_\infi$ is 
order-isomorphic to $\RR$ \cite[23.8]{Szm}. 
The corresponding ring epimorphism
\[
std:{}^\ast\RR_\fin\rTo\RR
\]
assigns to
every finite nonstandard real $x\in\RR_\fin$ its so-called standard part 
$std(x)\in\RR$.
We note also that $^\ast\RR$ contains canonically a copy of $\RR$
(represented by constant $K$-sequences of reals).
\[
 0\rTo\relax{}^*\RR_\infi\rTo\relax{}^*\RR_\fin \stackrel{\lTo}{\rTo}
 \RR\rTo 0
\]
\begin{Num}\textbf{Ultralimits and asymptotic cones\ }
\label{UltraLimitsAndAsymptoticCones}
For the remainder of this section, $D$ will denote a nonprincipal ultrafilter on the
countably infinite set $K$. We noted above that the ultraproduct
$\mathfrak{X}_D$ of a family of metric spaces is a generalized metric
space, whose distance function takes its values in the nonstandard reals%
\footnote{Generalized metric spaces, where the metrics take values in some
ordered abelian group
appear already in Hausdorff \cite[VIII \S5 III]{Haus}. 
See also \cite[Ch.~1.2]{Chi} \cite[II.1]{MorSha}.}.
Let $p$ be a point in the ultraproduct $X_D$ and let
$r$ be a positive nonstandard real. Let 
\[
X_D^{(p,r)}=\left\{q\in X_D\mid {\textstyle\frac1rd(p,q)\in{}^\ast\RR_\fin}\right\}.
\]
Then $\tilde d=std\circ \frac1rd$ is a real-valued pseudometric
\[
X_D^{(p,r)}\times
X_D^{(p,r)}\rTo^{\quad std\circ \frac1rd\quad }\RR,
\]
where points
with infinitesimal $\frac1rd$-distance have $\tilde d$-distance $0$.
Identifying points at $\tilde d$-distance $0$, 
we obtain a metric space which we
denote by $C(X_D,p,r)$,
\[
C(X_D,p,r)\times
C(X_D,p,r)\rTo^{\tilde d}\RR.
\]
This is the \emph{ultralimit}
of the family $(X_k)_{k\in K}$ with respect $p$ and $r$. 
If we represent the nonstandard real $r$ and the point
$p$ by $K$-sequences $(r_k)_{k\in K}$ and $(p_k)_{k\in K}$, respectively, then
these are the sequences of scaling factors and basepoints as in 
\cite[2.4.2]{KL} and \cite[7.5]{Roe}. In the notation of
these authors, we have
\[
C(X_D,p,r)=D\text-\lim\textstyle\frac1{r_k}X_k.
\]
If $r$ is infinite and all $X_k$ are
equal to one fixed metric space $X$, then $C(X_D,p,r)$ is the \emph{asymptotic cone} of 
$X$.
It is immediate that ultralimits commute with metric product decompositions.
From the surjection $std:{}^\ast\RR_\fin\rTo\RR$ we have in particular 
an isometry
\[
C({}^\ast\RR^n,p,r)\cong\RR^n
\]
for any choice of $p$ and $r$.
\end{Num}
Ultralimits have many good properties. For example, they are
spherically complete (and in particular complete, \cite[\S32]{Copson}).
Spherical completeness means that every nested sequence 
$\bar B_{t_0}(x_0)\supseteq \bar B_{t_1}(x_1)\supseteq \bar B_{t_2}(x_2)\supseteq 
\bar B_{t_3}(x_3)\supseteq \cdots$ of closed balls has nonempty intersection.
\begin{Lem}
The ultralimit $C(X_D,p,s)$ of any countably infinite family of metric spaces $(X_k)_{k\in K}$
is spherically complete.

\proof
Let $\bar B_{t_0}(x_0)\supseteq \bar B_{t_1}(x_1)\supseteq \bar B_{t_2}(x_2)\supseteq 
\bar B_{t_3}(x_3)\supseteq \cdots$ be a nested sequence
of closed balls in $C(X_D,p,r)$.
We consider the canonical surjection 
\[
X_D\supseteq X_D^{(p,r)}\rTo C(X_D,p,r)
\]
which identifies points whose $\frac1rd$-distance is infinitesimal.
For each $x_n$ we choose a preimage $z_n\in X_D^{(p,r)}$.
For $k,\ell\in\NN$ we put 
\[\textstyle
E_{k,\ell}=\left\{z\in X_D\mid \frac1rd(z,z_k)\leq t_k+\frac1{\ell+1}\right\}.
\]
This family has the finite intersection property.
Indeed, we have $z_m\in E_{k,\ell}$ for all $m\geq k$.
Adding constant symbols $z_n$ and $t_n$ to our
language $\mathcal L_\ms$, we may apply Theorem \ref{overspill}.
The countable set of formulas
\[\textstyle
\phi_{k,\ell}(x)=\biggl(\frac1rd(x,z_k)\leq t_k+\frac1{\ell+1}\biggr)
\]
can be satisfied by a single element $z\in X_D$, that is, 
$z\in \bigcap_{k,\ell}E_{k,\ell}$. 
The image of this element in $C(X_D,p,r)$ is contained in
$\bigcap_{n\geq 0}\bar B_{r_n}(x_n)$.
\qed
\end{Lem}
Now we suppose that each space $X_k$ is CAT$(\kappa_k)$.
Let $\kappa$ denote the nonstandard real corresponding to
the $K$-sequence $(\kappa_k)$. We also fix a positive
nonstandard real $r$. The next result is well-known
\cite[3.6]{KaLe} \cite[2.2.4]{KL}; from the viewpoint
of ultraproducts and \L os' Theorem, it is almost trivial.
\begin{Lem}
\label{UltralimitsPreserveCurvature}
For every real number $\lambda$ with $\lambda\geq \sqrt r\kappa$,
the space $C(X_D,p,r)$ is CAT$(\lambda)$. In particular,
$C(X_D,p,r)$ is CAT$(0)$ if each $X_k$ is CAT$(0)$.

\proof
This is obvious from \L os' Theorem \ref{LosThm} and the
fact that the CAT$(\kappa)$ condition can be stated in
the language $\mathcal L_\ms$. The details are as follows.
As in the standard definition of the CAT$(\kappa)$ condition
\cite[II.1]{BH}, 
we put $D_\kappa=\infty$ for $\kappa\leq 0$ and 
$D_\kappa=\frac\pi{\sqrt\kappa}$
for $\kappa>0$. For nonstandard reals $s<t$ we define the
`nonstandard interval' 
\[
[s,t]=\{x\in{}^\ast\RR\mid s\leq x\leq t\},
\]
and we define `nonstandard geodesics' in $X_D$ in the obvious way
as isometric injections of such nonstandard intervals.
It is clear from the fact that the $X_k$ are CAT$(\kappa_k)$ spaces
and from \L os' Theorem that
any two points $u,v\in X_D$ with $d(u,v)<D_\kappa$ can
be joined by a nonstandard geodesic.
We put $\eps=0$ if $\kappa=0$ and $\eps=\kappa/|\kappa|$ else.
Let $M^\eps$ denote the simply connected
complete Riemannian surface of constant sectional curvature
$\eps$, with its metric $d^\eps$, and let $M_D^\eps$ denote its
ultrapower. 
For $\kappa\neq 0$ we define a metric $d^\kappa$ on $M_D^\eps$ by
$d^\kappa=\frac1{\sqrt{|\kappa|}}d^\eps$.
Again, we see from \L os' Theorem \ref{LosThm}
that triangles of perimeter less than $2D_\kappa$ in $(X_D,d)$ are not thicker
than triangles in $(M^\eps_D,d^\kappa)$.
Thus $X_D$ is a `nonstandard CAT$(\kappa)$ space': it satisfies the
usual CAT$(\kappa)$ comparison triangle condition, except
that the model space $M^\eps$ is replaced by its nonstandard version $M^\eps_D$.
If we rescale the metric on $X_D$ by the factor $\frac1r$, we
obtain a nonstandard CAT$(\sqrt r\kappa)$ space.
If we identify in this space $(X_D,\frac1rd)$
points at infinitesimal distance, nonstandard geodesics
become ordinary geodesics. The CAT$(\lambda$) inequality also
remains valid for all real $\lambda\geq \sqrt r\kappa$.
Thus $C(X_D,p,r)$ is a CAT$(\lambda)$ space.
\qed
\end{Lem}
Asymptotic cones may be used to `make large-scale Lipschitz maps continuous' as follows.
\begin{Lem}
\label{CoarseToLipschitz}
Suppose that we have a countably infinite family of large-scale Lipschitz maps
\[
f_k:X_k\rTo Y_k,\quad\text{with}\quad d(f_k(u),f_k(v))\leq s_kd(u,v)+t_k.
\]
Let $s$ and $t$ denote the nonstandard reals represented
by the $K$-sequences $(s_k)$ and $(t_k)$ in $^*\RR$.
Assume that $s$ is finite with standard part $std(s)=s'\in\RR$.
Suppose that $r$ is a positive nonstandard real.
If $\frac tr$ is infinitesimal, then $f$ induces
an $s'$-Lipschitz map 
\[
C(f):C(X_D,p,r)\rTo C(Y_D,f_D(p),r).
\]
This holds in particular if the $s_k$ and $t_k$ are constant and
$r$ is any positive infinite nonstandard real.

\proof
We have $\frac1rd(f_D(u),f_D(v))\leq \frac1rd(u,v)+\frac tr$
by \L os' Theorem.
Identifying points at infinitesimal distance, we obtain a well-defined
Lipschitz map between the ultralimits.
\qed
\end{Lem}
\begin{Cor}\label{DimensionInvariance}
There is no quasi-isometric embedding $f:[0,\infty)\times\RR^n\rTo\RR^n$.

\proof
Otherwise, we could choose an infinite positive nonstandard real $r$ to
obtain a continuous injection $C(f)$
\begin{diagram}[height=1.5em]
C(([0,\infty)\times\RR^n)_D,0,r)&&\rInto^{C(f)}&& C(\RR^n_D,0,r)\\
\dTo^\cong&&&&\dTo^\cong\\
[0,\infty)\times\RR^n&&&& \RR^n.
\end{diagram}
Such a map would embed a closed $n+1$-cube homeomorphically into
$\RR^n$, which is impossible \cite[7.1.1,7.3.20]{Engelking}.
\qed
\end{Cor}

\section{Coarse equivalences of trees}
In this section we prove Theorem I from the introduction.
We first investigate the structure of group actions on
trees that are  $2$-transitive on the ends.
\begin{Def}
A metric space $T$ is called a \emph{tree} (or $\RR$-tree)
if it has the following two
properties.

\medskip\noindent
\textbf{(T1)} For any two points $x,y\in T$, there is a unique geodesic
$\gamma:[0,d(x,y)]\rTo T$ with $\gamma (0)=x$ and $\gamma(d(x,y))=y$.
We put $[x,y]=\gamma([0,d(x,y)])$.

\medskip\noindent
\textbf{(T2)} If $0<r<s$ and if $\gamma:[0,s]\rTo T$ is an injection such that
$\gamma|_{[0,r]}$ and $\gamma|_{[r,s]}$ are geodesics, then $\gamma$ is a geodesic.

\medskip\noindent
Trees are CAT$(\kappa)$ for every $\kappa\in\RR$, and in particular CAT$(0)$
\cite[II.1.15]{BH}. 
Basic references for trees are \cite{AlBa}, \cite{Chi} and \cite{MorSha}.
An $\RR$-tree with extensible geodesics is called a \emph{leafless tree}.
An \emph{apartment} in $T$ is an isometric image of $\RR$
(a geodesic line).
\end{Def}
If $z$ is in the geodesic segment $[x,y]$ but different from $x$ and $y$, 
we say that $z$ is \emph{between} $x$ and $y$. Given two geodesic segments
$[x,y]$ and $[x,y']$ in a tree, there is a unique
point $z$ with $[x,y]\cap[x,y']=[x,z]$ \cite[p.~30]{Chi}.
A point $z$ in the tree is called a \emph{branch point} if there are three points
$u,v,w$ distinct from $z$ such that $[u,v]\cap[v,w]\cap[w,u]=\{z\}$.
\begin{Num}\textbf{Ends\ }
A \emph{ray} in a tree $T$ is an isometric image of
$[0,\infty)$ (a geodesic ray). 
Two rays are called equivalent (or parallel) if their intersection is
again a ray; the resulting equivalence classes are called the
\emph{ends} of $T$. The set of all ends of $T$ is denoted by $\partial T$.
Given $x\in T$ and $u\in\partial T$, there is a unique geodesic
$\gamma:[0,\infty)\rTo T$ with $\gamma(0)=x$ whose image is in the
class of $u$ \cite[p.~60]{Chi}; we put $\gamma([0,\infty))=[x,u)$ and
$(x,u)=[x,u)-\{x\}$. 

Every apartment in $T$ determines two ends. Conversely,
if $u,v$ are distinct ends, then there is a unique apartment
whose ends are $u$ and $v$ \cite[p.~61]{Chi} and which we denote $(u,v)\subseteq T$.
If $A\subseteq T$ is an apartment and $z\in T$, then there exists a unique
point $\pi_A(z)$ in $A$ which has minimal distance from $z$ and
every geodesic segment $[z,x]$ with $x\in A$ contains $\pi_A(z)$
\cite[p.~61]{Chi}. This definition is compatible with \ref{CAT(k)Definition}.
\end{Num}
We now consider the following condition.

\medskip\noindent
{\boldmath\textbf{(2-$\partial$)}}
The group $G$ acts isometrically on the leafless metrically complete tree $T$, and
this action is $2$-transitive on $\partial T$: given ends
$u\neq v$ and $u'\neq v'$, there is an element
$g\in G$ with 
\[
g(u)=u'\quad\text{ and }\quad g(v)=v'.
\]
We note that such a group acts transitively on the set of apartments of $T$.
We now derive the structure of trees $T$ satisfying the condition (2-$\partial$).
We call a point $x$ in the tree \emph{$G$-isolated} if it is the 
unique fixed point of its stabilizer,
\[
T^{G_x}=\{x\}.
\]
For such a $G$-isolated point $x$, the stabilizer $G_x$ is by the Bruhat-Tits Fixed Point
Theorem \cite[II.2.8]{BH} a maximal bounded subgroup (i.e., $G_x$ is not properly contained in
a bounded subgroup of $G$). 
\begin{Prop}
\label{TransThm}
Let $u,v,w$ be three distinct ends of a tree $T$ satisfying (2-$\partial$)
and let $x$ be the branch point determined by these three
ends, $(u,v)\cap(u,w)=(u,x]$. Then there exists
an element $g\in G_x$ which fixes $u$ and maps $v$ to $w$.

\proof
Since $G$ is $2$-transitive on $\partial T$, we find
$h\in G_u$ such that $h(v)=w$.  
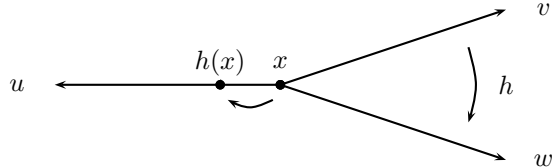
\begin{center}
\begin{pspicture}(6,2)
\psline[linestyle=solid]{<-*}(0,1)(3,1)
\psline[linestyle=solid]{<-*}(6,2)(3,1)
\psline[linestyle=solid]{<-*}(6,0)(3,1)
\psline[linestyle=solid]{*-*}(2.2,1)(2.2,1)
\pscurve[linestyle=solid]{->}(5.5,1.5)(5.6,1)(5.5,0.5)
\rput(-0.5,1){\small$u$}
\rput(6.5,2){\small$v$}
\rput(6.5,0){\small$w$}
\rput(6,1){\small$h$}
\rput(3,1.3){\small$x$}
\rput(2.2,1.3){\small$h(x)$}
\pscurve[linestyle=solid]{->}(2.9,0.8)(2.6,0.7)(2.3,0.8)
\end{pspicture}
\end{center}
Since $h(x)\in h(u,v)=(u,w)$, one of the
following three cases must hold:

(1) If $h(x)=x$ we are done, with $g=h$.

(2) If $h(x)\in(u,x)$, then $A=\bigcup\{h^{-n}(u,x]\mid n\in\NN\}$
is an $h$-stable apartment, on which $h$ induces a translation of length $d(x,h(x))$.
As all apartments are conjugate under $G$, we find an isometry $h'$ fixing $u$ and $w$
which maps $h(x)$ to $x$. Thus $g=h'h$ fixes $x$ and $u$ and maps $v$ to $w$.

(3) If $x\in(u,h(x))$, we put $A=\bigcup\{h^n(u,x]\mid n\in\NN\}$
and argue similarly as in (2).
\qed
\end{Prop}
\begin{Cor}
In a tree satisfying (2-$\partial$), all branch points are $G$-isolated.
\qed
\end{Cor}
We have the following a topological dichotomy.
\begin{Prop}
\label{DiscreteOrDense}
In a tree satisfying (2-$\partial$), the set of branch points is either closed and discrete, or it is
dense in every apartment.

\proof
Let $A=(u,v)$ be an apartment and assume that the set of branch points is not dense
in $A$. If $A$ contains no branch point, then $A=T$, and if $A$ contains
only one branch point, then every apartment contains only one branch point.
In either case, the claim follows.

Assume that $x\in A$ is a branch point and that $A$ contains another
branch point $y$. By \ref{TransThm} we find an
element $g\in G_x$ which interchanges $(u,x]$ and $(v,x]$.
Then $g(y)$ is also a branch point.
In this way, we get infinitely many branch points in $A$ at
uniform distance $d(x,y)$. Because the set of branch points in $A$ 
was assumed not to be dense, there has to be a minimal distance $t$ between these,
and they are distributed uniformly at this distance in $A$.

Since all apartments are conjugate, the set of branch points is either 
dense or discrete with uniform distance $t$ in every apartment.
In the discrete case, the minimal distance between two branch points is
$t$, so the set is closed and discrete in $T$.
\qed
\end{Prop}
\begin{Cor}
\label{TreeStructureCor}
Assume that the tree $T$ satisfies (2-$\partial$) and that the set of branch points is discrete. 
There are three possibilities for the structure of $T$.

\emph{Type (0)}. There are no branch points, $T\cong\RR$.

\emph{Type (I)}.
There is a single branch point. Then $T$ is the Euclidean cone
over its set of ends  $\partial T$, i.e. $T$ is a quotient of
$[0,\infty)\times\partial T$, where
$0\times v$ is identified with $0\times u$, for all $u,v\in\partial T$.

\emph{Type (II)}. There is an infinite discrete set of branch points.
Then $T$ is a simplicial metric tree, every vertex has valence at least $3$,
and all edges have the same length $t$.
\qed
\end{Cor}
This classification can be refined in terms of the $G$-action.
For type (0), $G$ induces a group $\{\pm 1\}\ltimes R$ on $T$,
where $R$ is a subgroup of $(\RR,+)$.
If $R\neq 0$, there are infinitely many $G$-isolated points.
For type (I), there is a unique $G$-isolated
point.
For type (II), $G_z$ acts transitively on the set of edges containing
$z$, for each branch point~$z$. Hence
there are two subcases: either $G$ acts transitively on the branch
points (vertices), or $T$ is bipartite and $G$ has two orbits on the
vertices. In the first case, the mid-points of the edges are 
$G$-isolated (and $G$ acts with inversion),
in the second case, only the branch points are $G$-isolated.
We note that for the types (I), (II), $T$ admits the structure of a
simplicial metric tree with a simplicial $G$-action.

We say that $T$ is of \emph{type (III)} if the set of branch points is
dense. By \ref{DiscreteOrDense}, the branch points are then dense in
every apartment.
\begin{Prop}
\label{AllPointsGIsolated}
Assume that the tree $T$ satisfies (2-$\partial$) and that the set of branch points is dense.
Then every point $x\in T$ is $G$-isolated.

\proof
Let $x\in T$. If $G_x$ fixes another point
$y\neq x$, then $G_x$ fixes the geodesic segment $[x,y]$.
There is a branch point $z$ between $x$ and $y$. By \ref{TransThm}
$y$ is not a fixed point of $G_{x,z}$. But $G_{x,y}=G_x\supseteq G_{x,z}$, a contradiction.
\qed
\end{Prop}
Every tree satisfying (2-$\partial$) corresponds
to one of the four types (0)-(III). Now we clarify the role of the maximal
bounded subgroups. We noted already that by the Bruhat-Tits Fixed Point
Theorem, the stabilizers of $G$-isolated points are maximal bounded subgroups.
\begin{Prop}
\label{MaxBoundedSubgroups}
Assume (2-$\partial$) and that $P\subseteq G$ is a maximal bounded subgroup.
Then $P$ is the stabilizer of a $G$-isolated point.

\proof
Assume this is false, so the fixed point set $T^P$ contains a geodesic segment
$[x,y]$ with $x\neq y$. If $[x,y]$ contains a branch point
$z$, then $G_z=P$ (by maximality), whence $T^P=T^{G_z}=\{z\}$, a contradiction.
From \ref{DiscreteOrDense} we see that the set of branch points 
is closed and discrete in $T$. By \ref{TreeStructureCor} we can assume that
$T$ is a simplicial tree. Then $P$ fixes some simplicial edge of $T$ elementwise.
Therefore it fixes also some branch point $z$, which is a contradiction
to \ref{TransThm}.
\qed
\end{Prop}
In a tree $T$ satisfying (2-$\partial$), 
let $i_G(T)$ denote the set of $G$-isolated points of $T$.
By the results above, the  set $i_G(T)$ corresponds bijectively to the set of maximal
bounded subgroups of $G$. With respect to the conjugation action of $G$ on subgroups, 
this correspondence is $G$-equivariant. Our aim is to show that $T$ can be recovered
from the $G$-actions on $i_G(T)$ and $\partial T$. We let 
\[
b(T)\subseteq i_G(T)
\]
denote the set of branch points and consider the different types (0)-(III) of trees.
By the \emph{combinatorial structure} of a tree we mean the underlying set $T$
together with the collection of all apartments in $T$ (without any metric).
First, we dispose of the two degenerate types (0) and (I).
\begin{Num}
\label{Recover_a_b}
Type (0) is characterized by $\#\partial T=2$ (and $b(T)=\emptyset$).
The tree $T\cong\RR$ is unique
up to isometry and the group induced by $G$ splits as a semidirect product
$\ZZ/2\ltimes R$, where $R$ is a subgroup of
$(\RR,+)$. Since $\RR$ contains subgroups which are
abstractly, but not topologically, isomorphic, the $G$-action on $T$ can
in general not be recovered from the $G$-action on $i_G(T)$ and $\partial T$.

Type (I) is characterized by $\#\partial T\geq 3$ and $\#i_G(T)=1$.
This determines both the combinatorial structure of the tree $T$ and the
$G$-action, but not the metric.
\end{Num}
In the remaining cases, both $\partial T$ and $i_G(T)$ are infinite.
This situation is much more rigid.
\begin{Num}
\label{Recover_c}
Suppose that $T$ is of type (II). Then $x\in i_G(T)$ is a branch point
if and only if $G_x$ has no orbit of length $2$ in $i_G(T)$.
So we can recover the set $b(T)$ of branch points in $i_G(T)$.
By \ref{TransThm}, two branch points $x,y$ are adjacent if and only if
the only branch points fixed by $G_{x,y}$ are $x$ and $y$, hence the simplicial
structure of $T$ can be recovered from the $G$-action on $i_G(T)$. 
We note that then $G_{x,y}$ has at most $3$ fixed points in
$i_G(T)$. Since all edges of $T$ have the same length, $T$ is determined as a
metric space up to a scaling factor.
\end{Num}
The following result shows how for type (III) apartments can be described
in terms of the $G$-action. Recall that $\pi_A:T\rTo A$ denotes the
retraction that sends $x$ to the nearest point $\pi_A(x)\in A$.
\begin{Lem}
\label{RecoverApts}
Assume that $T$ is of type (III). Let $A=(u,v)$ be an apartment.
Then $x$ is contained in $A$ if and only if $x$ is the unique fixed point of
$\bra{G_{x,u}\cup G_{x,v}}$.

\proof
Suppose that $x\not\in A$.
Then both $G_{x,u}$ and $G_{x,v}$ fix $\pi_A(x)\neq x$.
Now assume that $x\in A$ and that $z\neq x$.
If $\pi_A(z)=x$ then $x$ is a branch point and $G_{u,x}$ moves
$z$ by \ref{TransThm}. If $\pi_A(z)$ is between $x$ and $u$,
then there is a branch point between $x$ and $\pi_A(z)$ and
therefore $G_{v,x}$ moves $z$ by \ref{TransThm}. Similarly, $G_{u,x}$
moves $z$ if $\pi_A(z)$ is between $x$ and $v$.
\qed
\end{Lem}
\begin{Num}
\label{Recover_d}
If $T$ is of type (III), then $i_G(T)=T$ by \ref{AllPointsGIsolated}. 
The stabilizer $G_{x,y}$
of two distinct points $x,y$ fixes the infinite set $[x,y]\subseteq i_G(T)$, so
type (III) can be distinguished from type (II). By \ref{RecoverApts} we
can recover the apartments in $T$ from the group action.
Let $A=(u,v)$ be an apartment and let $z\in A$. Then 
$(u,z]=A\cap T^{G_{u,z}}$, so we can also recover the rays in $A$.
Therefore $G$ determines the topology of $A$.
Since the branch points are dense $A$, the group $G_{u,v}$
induces a group $H$ of translations on $A$ which has a dense
orbit. Up to a scaling factor, there is just one $H$-invariant
metric on $A$ which satisfies $d(h^n(x),y)=nd(x,y)$ for
all $h\in H$ and all $n\in\NN$.
So the metric on $A$ is determined up to scaling. Since all
apartments are conjugate, the metric on $T$ is unique up to
scaling.
\end{Num}
Summarizing these facts we have the following result.
\begin{Prop}
\label{GroupEncodesTree}
Assume that $T$ is a tree with $\#\partial T\geq 3$ and that $G$ acts on $T$
in such a way that (2-$\partial$) holds. 
Given the $G$-actions on $\partial T$ and $i_G(T)$, the combinatorial structure of the
tree $T$ is uniquely determined. If $T$ is of type
(II) or type (III), the metric is determined
up to a scaling factor.
\qed
\end{Prop}
Now we get to our main rigidity result for trees. Suppose that
\[
 f:T_1\rTo T_2
\]
is a coarse equivalence between two trees. We recall
a result about coarse equivalences between $\delta$-hyperbolic spaces
which goes back to M.~Morse. It says that the coarse
image of a geodesic is Hausdorff close to a geodesic
\cite[III.H.1.7]{BH} \cite[1.3.2]{BS}. Since trees
are $\delta$-hyperbolic for every $\delta>0$, it follows that there is a positive constant $r_f>0$
such that the image $f(A)$ of an apartment $A\subseteq T_1$
has Hausdorff distance at most $r_f$ from a unique apartment $A'\subseteq T_2$.
We put $f_*A=A'$. Note that $f_*$ just maps apartments to apartments, it is
not a map between the trees. 
\begin{Num}\label{Quasigeodesics2}
For our application to Euclidean buildings we have to consider the more
general situation of a coarse equivalence 
\[
f:T_1\times\RR^{m_1}\rTo T_2\times\RR^{m_2}.
\]
A leafless $\RR$-tree is
a $1$-dimensional Euclidean building, so we may apply the
higher dimensional Morse Lemma \ref{HigherDimensionalMorseLemma}.
It follows that there is a constant
$r_f>0$ such that for every apartment $A\subseteq T_1$
there is a unique apartment $A'\subseteq T_2$ such that
$f(A\times\RR^{m_1})$ has Hausdorff distance at most $r_f$ from
$A'\times\RR^{m_2}$. Moreover, 
\[
m_1=m_2
\]
by \ref{CoarseEqivalencesPreserveEuclideanFactors}.
We define a map $f_*$ from the set of apartments in $T_1$ to the set of apartments
in $T_2$ by putting
\[
f_*A=A'.
\]
\end{Num}
If $g$ is a coarse inverse for $f$, then
$g_*$ is an inverse for $f_*$.
We also need the following auxiliary result on trees.
\begin{Lem}
\label{BoundedNbhd}
Let $\F $ be a collection of apartments in a tree $T$ and let
$r>0$. If the subset $X=\bigcap\{B_r(A)\mid A\in\F \}\subseteq T$
is nonempty and unbounded, then
the apartments in $\F $ have a common end.

\proof
The result is a special case of \ref{AptCplxHom} below, but we give
a direct proof. Let $(u,v)\in\F$ and choose a sequence
$x_n\in X$ such that $\pi_{(u,v)}(x_n)$ converges to one end of
$(u,v)$, say $u$. This is possible since $X\subseteq B_r((u,v))$ is
unbounded.
Let $A\in\F$. Then $d(\pi_A(x_n),\pi_{(u,v)}(x_n))\leq 2r$ and the
unbounded sequence $\pi_A(x_n)$ subconverges to an end $w$ of $A$.
If $w\neq u$, then the set $\{d(\pi_A(x_n),\pi_{(u,v)}(x_n))\mid n\in\NN\}$ would
be unbounded. Thus $w=u$.
\qed
\end{Lem}
In the previous lemma the
end is unique, unless $\F $ consists of a single apartment. 
\begin{Prop}
\label{BoundaryLemma}
If $T_1$ has at least one branch point and if $f:T_1\times\RR^m\rTo T_2\times\RR^m$
is a coarse equivalence, then the map $f_*$ between the sets
of apartments of the trees
induces a bijection $f_*:\partial T_1\rTo\partial T_2$
between the ends of the trees, in such a way that
$f_*(u,v)=(f_*u,f_*v)$.

\proof
Let $r_f>0$ be as in \ref{Quasigeodesics2} and let
$\mathcal F$ be a finite collection of apartments in $T_1$ having an
end $u$ in common. Put $(u,x]=\bigcap\{A\mid A\in\mathcal F\}$ and 
$Y=\bigcap \{B_{r_f}(f_*A)\mid A\in\mathcal F\}$.
Then $f((u,x]\times\RR^m)\subseteq Y\times\RR^m$.
If the set $f_*\mathcal F$ has no common end, then
$Y$ is bounded by \ref{BoundedNbhd},
so $Y\times\RR^m$ is quasi-isometric to $\RR^m$.
We obtain then a quasi-isometric embedding of $[0,\infty)\times\RR^m$ into
$\RR^m$, which is impossible by \ref{DimensionInvariance}.
Therefore $Y$ is unbounded and 
$f_*\mathcal F$ has a common end. 
Applying the same argument to a coarse inverse $g$ of $f$, we see that
a finite set $\mathcal F$ has an end in common if and only if 
$f_*\mathcal F$ has an end in common.

For $i=1,2$, consider the simplicial complex $AC(T_i)$
whose simplices are the
finite collections of apartments in $T_1$ having a common end.
The ends of $T_i$ correspond to the maximal complete subcomplexes
of $AC(T_i)$ (a simplicial complex is called complete if any two
simplices are contained in some simplex).
We have shown that $f_*$ induces an isomorphism $f_*:AC(T_1)\rTo AC(T_2)$.
It follows that $f_*$ induces a bijection $\partial T_1\rTo\partial T_2$.
\qed
\end{Prop}
The complex $AC(T_i)$ is a special case of the apartment complex of a
spherical building;  see \ref{AptCplxDef}. It is the
nerve of the covering of $\partial T_i$ given by all two-element subsets.
The last paragraph of the proof is a special case 
of~\ref{AptCplxProp} below.
The next result is a special case of \cite[8.3.11]{KL}.
\begin{Prop}
\label{BoundedThm}
Let $G$ be a group acting isometrically on two metrically complete leafless
trees $T_1,T_2$.
Assume that $f:T_1\times\RR^m\rTo T_2\times\RR^m$ is a coarse equivalence and that
the induced map $f_*$ on the apartments is $G$-equivariant.
If $T_1$ has at least three ends, then a subgroup $P\subseteq G$
has a bounded orbit in $T_1$ if and only if it has a bounded orbit
in $T_2$.

\proof
Suppose that $P\subseteq G$ has a bounded orbit in $T_1$. Let $x$ be a branch point in $T_1$ and
consider the bounded orbit $P(x)$. We put
\[
f(x\times0)=y\times q
\]
and we show that $y$ has a bounded orbit $P(y)$ in $T_2$.

Let $\F $ denote the set of all apartments in $T_1$ which
intersect the orbit
$P(x)$ nontrivially. This set $\F$ is obviously $P$-invariant
and has no common end (because $x$ is a branch point).
Let $\F'=f_*\F $ denote the corresponding set of
apartments in $T_2$. Since we assume that $f_*$ is $G$-equivariant,
$\F'$ is also $P$-invariant, and the apartments in $\F'$
have no common end by \ref{BoundaryLemma}.
For $s>0$ consider the $P$-invariant set
\[
X_s=\bigcap\{B_s(A')\mid A'\in \F'\}\subseteq T_2.
\]
By \ref{BoundedNbhd}, the set $X_s$ is bounded (or empty).
Let $r'=\sup\{d(x,p(x))\mid p\in P\}$. For every $A\in\F$ we have $d(x,\pi_A(x))<r'$.
Put $f(\pi_A(x)\times0)=y'\times p'$.
Then $d(y',\pi_{f_*A}(y'))\leq r_f$ by \ref{Quasigeodesics2}.
If $\rho$ denotes the control function for $f$, then $d(y,y')\leq\rho(r')$. This holds for all
$A\in\F$, whence 
\[
y\in X_{\rho(r')+r_f}.
\]
As this set is bounded and $P$-invariant, $P(y)$ is bounded.
If $g$ is a coarse inverse of $f$, then $g_*$ is $G$-equivariant
(because it is the inverse of the equivariant map $f_*$),
so we obtain the converse implication by the same arguments.
\qed
\end{Prop}
Now we prove our main result on trees, which implies Theorem~I in the
introduction.
The map $f_*:\partial T_1\rTo\partial T_2$
is defined as in \ref{BoundaryLemma}.
\begin{Thm}
\label{MainTreeProp}
Let $G$ be a group acting isometrically on two
metrically complete leafless
trees $T_1,T_2$, with $\#\partial T_1\geq 3$ and
assume that the action of $G$ on $\partial T_1$
is $2$-transitive. Suppose that there is a coarse equivalence
\[
f:T_1\times\RR^m\rTo T_2\times\RR^m
\]
and that the induced map
$f_*:\partial T_1\rTo\partial T_2$ is $G$-equivariant.
Then we have the following.

(i) After rescaling the metric on $T_2$ by a constant $r>0$, there is
a $G$-equivariant isometry 
\[
\bar f:T_1\rTo T_2.
\]
For every apartment $A\subseteq T_1$ we have $\bar f(A)=f_*A$.  
If $T_1$ has at least $2$ branch points, then both $\bar f$ and $r$ are unique.

(ii) Put $f(x\times p)=f_1(x\times p)\times f_2(x\times p)$.
If $T_1$ has at least $2$ branch points, then there is
a constant $s>0$ such that $d(f_1(x\times p),\bar f(x))\leq s$
holds for all $x\times p\in T_1\times\RR^m$.
The constant $s$ depends only on $T_1$, the control function $\rho$, $a$ and the
constant $r_f$ from \ref{Quasigeodesics2}.
In particular, $f$ and $\bar f$ have finite distance
if $m=0$.

(iii) If $f$ is a rough isometry (see \ref{QIDef}), we may put $r=1$.

\proof
(i)
By \ref{BoundedThm} both $G$-actions have the same set of
maximal bounded subgroups. These subgroups correspond by
\ref{MaxBoundedSubgroups} to the $G$-isolated points.
Therefore we have an equivariant
bijection
\[
\bar f:i_G(T_1)\rTo i_G(T_2).
\]
The types
(I), (II), and (III) can be distinguished by the $G$-action on
$i_G(T_1)$; see \ref{Recover_a_b}, \ref{Recover_c} and \ref{Recover_d}
(our assumptions exclude trees of type (0)).
The combinatorial structure is also encoded in the
$G$-action, as we  noted in \ref{GroupEncodesTree}. 
Also, we can rescale the metric on $T_2$ by a constant $r>0$
in such a way that $\bar f:i_G(T_1)\rTo i_G(T_2)$
extends to an equivariant isometry $T_1\rTo T_2$ which we also denote
by $\bar f$. From the construction it is clear that
$\bar f(A)=f_*A$. For trees of type (II) and (III), $\bar f$ and
$r$ are unique by~\ref{GroupEncodesTree}.

(ii)
Let $z\in b(T_1) $ and
put $r_1=1+\inf\{d(x,y)\mid x,y\in b(T_1),x\neq y\}$. Then
$T_1$ is covered by the $G$-translates of $B_{r_1}(z)$.
Let $\F $ be the collection of all apartments of $T_1$ containing $z$.
By \ref{Quasigeodesics2} there is a constant $r_f>0$ such that
$f_1(z\times p)\in \bigcap\{B_{r_f}(f_*A)\mid A\in\F \}=
\bigcap\{B_{r_f}(\bar f(A))\mid A\in\F \}=B_{r_f}(\bar f(z))$.
It follows that $d(f_1(x\times p),\bar f(x))\leq r r_1+r_f+\rho(r_1)$
for general $x\in T_1$.

(iii)  For trees of type (I) it
is clear that no rescaling is necessary in order to find
$\bar f$. 
Suppose that $f$ is a rough isometry, with control function
$\rho(t)=t+b$. For trees of type (II) and (III) we have by (ii)
$d(\bar f(x),\bar f(x'))\leq d(x,x')+b+2(rr_1+r_f+b)$.
On the other hand, $d(\bar f(x),\bar f(x'))=rd(x,x')$.
Since $T_1$ is unbounded, $r\leq 1$.
Applying the same argument to $\bar f^{-1}$, we see that $r=1$.
\qed
\end{Thm}
For the special case of a coarse equivalence between locally finite 
simplicial trees, a similar result is proved in \cite[4.3.1]{Le}.

\section{Spherical buildings}
In this section we record some basic notions and facts about spherical
buildings. Everything we need can be found in 
\cite{Bro}, \cite{Ron}, \cite{TiLNM} and \cite{WeSph}.
For our present purposes it is convenient to view buildings as
simplicial complexes. This is essentially Tits' approach 
in \cite{TiComo}; see also \cite{CLT}.
\begin{Num}\textbf{Simplicial complexes\ }
\label{SimplicialComplexes}
Let $V$ be a set and $S$ a collection of finite subsets
of $V$. If $\bigcup S=V$ and if $S$ is closed under going down
(i.e. $a\subseteq b\in S$ implies $a\in S$), then the poset
$(S,\subseteq)$ is called a \emph{simplicial complex}.
The $k+1$-element subsets $a\in S$ are called $k$-\emph{simplices}.
More generally, any poset isomorphic to such a poset $(S,\subseteq)$
will be called a simplicial complex.
The \emph{join} $S*T$ of two simplicial complexes $S,T$ is
the product poset; it is again a simplicial complex.
Homomorphisms between simplicial complexes are defined in the
obvious way as order-preserving maps which do not raise the dimension of
simplices \cite[7A.1]{BH}. A homomorphism which maps $k$-simplices to
$k$-simplices is called \emph{non-degenerate}. 
The
\emph{geometric realization} $|S|$ of the simplicial complex $S$
is the set
\[
|S|=
\left\{
p:V\rTo\relax[0,1]\mid p^{-1}(0,1]\in S\text{ and }{\textstyle\sum_{v\in V}p(v)=1}
\right\}
\]
Endowed with the weak topology, $|S|$ is then a CW complex.
Moreover, 
\[
|S*T|\cong|S|*|T|,
\]
where the right-hand side
is the topological join.
In the case of spherical buildings, $|S|$ is often endowed with
a stronger, metric topology \cite[II.10A]{BH}. A result due to Dowker says
that the identity map is a homotopy equivalence between these two topologies
\cite[I.7]{BH}.
\end{Num}
\begin{Num}\textbf{Coxeter groups and buildings\ }
Let $(W,I)$ be a Coxeter system \cite[IV.3]{Bourbaki} \cite[II.4]{Bro}. 
Thus $W$ is a group with a (finite)
generating set $I$ consisting of involutions and a presentation
of the form 
\[
W=\left\langle I\ \left|\ (ij)^{ord(ij)}=1\text{ for all }i,j\in I\right\rangle\right..
\]
For a subset 
$J\subseteq I$ we put $W_J=\bra{J}$. If $J$ is nonempty, then $(W_J,J)$ is again a Coxeter
system \cite[IV.8]{Bourbaki}. The poset 
\[
\Sigma=\Sigma(W,I)=\bigcup\{W/W_J\mid J\subseteq I\},
\]
ordered by
reversed inclusion, is a simplicial complex, the \emph{Coxeter complex}
\cite[III.1]{Bro}. The \emph{type} of a simplex $wW_J$ is
$t(wW_J)=I- J$. The type function may be viewed as a non-degenerate
simplicial epimorphism from $\Sigma$ to the power set $2^I$ of $I$
(viewed as a simplicial complex).

A \emph{building} $B$ is a simplicial complex together with
a collection $Apt(B)$ of subcomplexes, called \emph{apartments}, which are
isomorphic to a fixed Coxeter complex $\Sigma$. The apartments
have to satisfy the following two compatibility conditions.

\medskip\noindent
\textbf{(B1)}
For any two simplices $a,b\in B$, there is an apartment $A\in Apt(B)$
containing $a,b$.

\medskip\noindent
\textbf{(B2)}
If $A,A'\subseteq B$ are apartments containing the simplices $a,b$, then there is
a (type preserving) simplicial isomorphism $A\rTo A'$ fixing $a$ and $b$.

\medskip\noindent
The type functions of the apartments are
therefore compatible and extend to a non-degenerate surjective simplicial map
$t:B\rTo 2^I$, the type function of the building.
The cardinality of $I$ is the \emph{rank} of the building,
\[
\#I=\mathrm{rank}(B)=\dim(B)+1.
\]
A building of rank $1$ is just a set (of cardinality at least $2$),
the apartments are the two-element subsets.

The maximal simplices in a building are called \emph{chambers}.
We denote by $Cham(B)$ the set of all chambers.
Every simplex
of a building is contained in some chamber (so buildings are
\emph{pure} simplicial complexes). Recall that the \emph{dual graph} of a
pure complex is the graph whose vertices are the maximal
simplices and whose edges are the simplices of codimension $1$;
this is the \emph{chamber graph} of $\Delta$.
A \emph{gallery} is a simplicial path in the chamber graph, and a
nonstammering gallery is a path where consecutive chambers are
distinct. The chamber graph of a building is always
connected. A \emph{minimal gallery} is a shortest path in the
chamber graph.
\end{Num}
A building is called \emph{thick}
if every non-maximal simplex is contained in at least $3$
distinct chambers (it is always contained in at least $2$ distinct
chambers). If every simplex of codimension $1$ is contained in
exactly two chambers, the building is called \emph{thin}.
Thin buildings are Coxeter complexes.
We allow non-thick buildings (in \cite{TiLNM}, all
buildings are assumed to be thick, but the results from \cite{TiLNM}
which we collect in this section hold for non-thick buildings as well).
\begin{Num}\textbf{Residues and panels\ }
Let $a\in B$ be simplex of type $J$. The \emph{residue} of $a$ is the poset
\[
\Res(a)=\{b\in B\mid b\supseteq a\}.
\]
If $a$ is not a chamber, then this poset is again a
building, whose Coxeter complex is modeled on $W_{I- J}$
\cite[3.12]{TiLNM}.
If $a$ is a simplex of codimension $1$ and type $I-\{j\}$,
then $\Res(a)$ is called a \emph{$j$-panel}.

There is an order-reversing poset isomorphism between the simplicial complex
$B$ and the set of all residues in $B$.
Residues can also be defined in terms of the chamber graph, viewed
as an edge-colored graph. Basically, this is a dictionary which allows
the passage from buildings, viewed as simplicial complexes (as in \cite{TiLNM})
to buildings viewed as edge colored graphs (or chamber systems)
(as in \cite{WeSph}). In view of this correspondence, we call a simplex of type
$I-\{j\}$ also a $j$-panel.
\end{Num}
The join of two buildings
is again a building. Conversely, a building decomposes as a join
if its Coxeter group is decomposable (i.e. if $I\subseteq W$
decomposes into two subsets which centralize each other).
\begin{Num}\textbf{Spherical buildings\ }
A Coxeter complex $\Sigma$ is called \emph{spherical} if it is
finite. Then the geometric realization $|\Sigma|$ is a combinatorial
sphere of dimension $\#I-1$.
A building is called \emph{spherical} if its apartments are finite
(the building may nevertheless be infinite).
Two simplices $a,b$ in a spherical Coxeter complex $\Sigma$ are called
\emph{opposite} if they are interchanged by the antipodal map
(the opposition involution) of the sphere $|\Sigma|$ \cite[3.22]{TiLNM}.
In a spherical building, two simplices are called opposite if they
are opposite in some (hence every) apartment containing them.
\end{Num}
\begin{Num}\textbf{Thick reductions\ }
\label{ThickRedSph}
A non-thick
spherical building $B$ can always be `reduced' to a thick building as
follows. There exists a thick spherical building $B_0$ such that
$B$ is a simplicial refinement of a join
$\SS^0*\cdots *\SS^0*B_0$ (we view $\SS^0$ as a thin spherical building
of rank $1$) \cite{Cap} \cite{CL} \cite[3.7]{KL} \cite[3.8]{KramerCR}
\cite{Scha}.
For the geometric realization, we have then
\[
|B|=\SS^k*|B_0|,
\]
where $k$ is the number of $\SS^0$-factors in the join.
So non-thick spherical buildings are suspensions of thick spherical buildings.
The geometric realization of a thin spherical building is a sphere.
\end{Num}
The following lemma will be used later.
\begin{Lem}
\label{ThickApartmentsLemma}
Let $B$ be a spherical building and let $A\subseteq B$
be an apartment. If every panel $a\in A$ is contained in
at least three different chambers of $B$, then $B$ is
thick.

\proof
Let $a$ be an arbitrary panel in $B$, and let
$c_0,\cdots,c_k$ be a shortest gallery with the property that
the first chamber $c_0$ contains $a$ and the last chamber
$c_k$ is in $A$. This gallery can be continued inside $A$
as a minimal gallery until it reaches a panel $b\in A$ which is opposite $a$.
Then there
is a bijection $[a;b]$ between $\Res(a)$ and $\Res(b)$; see \ref{Perspectivities} below.
It follows that $a$ is contained in at least three different chambers.
\qed
\end{Lem}
We will see that Euclidean buildings give rise to a family of building epimorphisms. 
The following fact is useful; see \cite[\ref{MaxBoundedSubgroups}]{BK}.
\begin{Lem}
\label{LocalSurjectivity}
Let $B,B'$ be spherical buildings of the same type and let
$\phi:B\rTo B'$ be a type-preserving simplicial map.
Then $\phi$ is 
surjective if and only if its restriction to every
panel is surjective.

\proof
Since the chamber graph of $B'$ is connected, it is clear that the
local surjectivity condition implies global surjectivity.
Conversely, suppose that $\phi$ is surjective. Let
$a',b'\in B'$ be $i$-adjacent chambers and let $a$ be a
preimage of $a'$. We have to find a preimage $b$ of
$b'$ which is $i$-adjacent to $a$. Let $c'$ be opposite
$a'$ and let 
\[
a'\rLine_i b'\rLine\cdots\rLine c'
\]
be
a minimal gallery which we denote by $\gamma'$.
Let $c$ be a preimage of $c'$. Since
$\phi$ does not increase distances in the chamber graph, $c$
is opposite $a$. Thus there is a gallery 
\[
a\rLine_i b\rLine\cdots\rLine c
\]
in $B$ of the same type as $\gamma'$,
which we denote by $\gamma$. Since $\gamma'$ is the unique gallery
of its type in $B'$ from $a'$ to $c'$, it follows that
$\phi(\gamma)=\gamma'$. In particular, $\phi(b)=b'$.
\qed
\end{Lem}
\begin{Cor}
\label{ThicknessPullsBack}
Let $B,B'$ be spherical buildings of the same type and let
$\phi:B\rTo B'$ be a type-preserving simplicial surjective map.
If $B'$ is thick, then $B$ is also thick.
\qed
\end{Cor}
A thick spherical building is determined by its apartment
complex which we introduce now. This complex appeared already
in \ref{BoundaryLemma} for the special case of rank $1$ buildings. 
We will see later that the apartment complex of the spherical
building at infinity is a coarse invariant of a Euclidean building.
\begin{Num}\textbf{The apartment complex\ }
\label{AptCplxDef}
Let $B$ be a building and $Apt(B)$ its set of apartments. The 
\emph{apartment complex} $AC(B)$ is the simplicial complex whose
simplices are finite subsets $A_1,\ldots,A_k$ of apartments,
with $A_1\cap\cdots\cap A_k\neq\emptyset$ (in other words,
$AC(B)$ is the \emph{nerve} of the covering $Apt(B)$ of $B$).
\end{Num}
If $B$ is thick, then every simplex $a\in B$ can be written as
an intersection of finitely many apartments. 
\begin{Prop}
\label{AptCplxProp}
Let $B_1,B_2$ be thick spherical buildings and let
$\phi:AC(B_1)\rTo AC(B_2)$ be a simplicial isomorphism.
Then there is a unique simplicial isomorphism $\Phi:B_1\rTo B_2$ such that
$\phi(A)=\Phi(A)$ for all $A\in Apt(B_1)$.

\proof
The map $\phi$ sends sets of apartments with the finite intersection
property to sets with the finite intersection property. 
Because $B_1$ is thick and has finite apartments,
the maximal subsets with
the finite intersection property in $Apt(B_1)$ are precisely the sets
$\mathcal S_v=\{A\in Apt(B_1)\mid v\in A\}$ for all vertices $v$ of $B$.
Therefore $\phi$ induces a bijection $\Phi$ between the vertices of
the buildings. If $v\in A$, then $\Phi(v)\in\phi(A)$. Since
$B_1$ is thick, two vertices $u,v$ in $B_1$ are adjacent if and
only if the following holds: no other vertex $w$ is in the
intersection of all apartments containing $u$ and $v$. This can be expressed
in $AC(B_1)$ as follows: if $\mathcal S_w\supseteq \mathcal S_u\cap
\mathcal S_v$, then $\mathcal S_w=\mathcal S_u$ or
$\mathcal S_w=\mathcal S_v$. So $\Phi$ preserves the $1$-skeleton of
$B_1$. But every building is the flag complex of its $1$-skeleton
\cite[3.16]{TiLNM},
therefore $\Phi$ is simplicial.
\qed
\end{Prop}
In the previous proof, thickness is essential. It is clear that
the thick reduction (see \ref{ThickRedSph}) of a spherical
building has the
same apartment complex as the building itself.
We remark that a simplicial isomorphism between two spherical
buildings maps apartments to apartments, even if it is not type-preserving.
\begin{Num}\textbf{Projections in buildings\ }
\label{Projections}
Let $c$ be a chamber and $a$ a simplex in a building.
Then there is a
unique chamber $d$ in $\Res(a)$ which has minimal distance
from $c$ (with respect to the distance in the chamber graph),
and which is denoted $d=\proj_ac$. If $b$ is a simplex
then $\proj_ab$ is defined to be the simplicial intersection
of the chambers $\proj_ac$, where $c$ ranges over all
chambers containing $b$ \cite[3.19]{TiLNM}.
\end{Num}
If $a,b$ are opposite simplices in a spherical building,
then $\proj_a:\Res(b)\rTo \Res(a)$ is a
simplicial isomorphism \cite[3.28]{TiLNM}.
The following observations are due to Knarr \cite{Kn}
and Tits \cite{TiComo}. They were rediscovered by 
Leeb \cite[Ch.~3]{Le}.
\begin{Num}\textbf{Perspectivities\ }
\label{Perspectivities}
Let $B$ be a spherical building. We noted already that if $a,b$ are
opposite simplices in $B$, then $\proj_b:\Res(a)\rTo \Res(b)$ is a
building isomorphism
(not necessarily type preserving) between $\Res(a)$ and $\Res(b)$.
We denote this isomorphism by 
\[
[b;a]:\Res(a)\rTo \Res(b)
\]
and
call it a \emph{perspectivity}. A concatenation of perspectivities is called a 
\emph{projectivity}; we write 
\[
[c;b]\circ[b;a]=[c;b;a]:\Res(a)\rTo \Res(c)
\]
etc.
The inverse of $[b;a]$ is $[a;b]$. A projectivity is called \emph{even}
if it can be written as a composition of an even number of perspectivities.
\end{Num}
\begin{Num}\textbf{Projectivities\ }
\label{Projectivities}
Recall that a \emph{groupoid} is a small category where every arrow is an
isomorphism.
The \emph{projectivity groupoid} $\Pi_B$ of a spherical building $B$ is
the category whose objects are the simplices of $B$, and whose morphisms
are projectivities. It is closely related to the \emph{opposition graph}
$\Opp(B)$ whose vertices are the simplices of $B$ and whose edges
are unordered pairs of opposite simplices. Every simplicial path in
$\Opp(B)$ induces a projectivity.
We denote by 
\[
\Pi_B(a)=\mathrm{Hom}_{\Pi_B}(a,a)
\]
the group of all
automorphisms of $\Res(a)$ induced by maps in $\Pi_B$.
The subgroup $\Pi_B(a)^+$ consisting of all even projectivities
is a normal subgroup of index $1$ or $2$ in $\Pi_B(a)$.
If $f:B_1\rTo B_2$ is an isomorphism of spherical buildings, then
$f$ induces an isomorphism between $\Pi_{B_1}$ and $\Pi_{B_2}$
in the obvious way.
\end{Num}
The following result is essentially Knarr's~\cite[1.2]{Kn}.
We use several facts about galleries and distances
which can all be found in \cite{WeSph}.
\begin{Thm}
\label{KnarrThm}
Suppose that $B$ is a thick spherical building and that $r$ is an
$i$-panel.
If $i$ is not an isolated node in the Coxeter diagram of $B$,
then $\Pi_B(r)^+$ is a $2$-transitive permutation group on $R=\Res(r)$.

\proof Let $a,b,b'$ be three distinct chambers in $R$. We construct a projectivity
which fixes $a$ and maps $b$ to $b'$. Let 
$j$ denote the type of a neighboring node of $i$, i.e.
$ij\neq ji$ in $W$. We choose a nonstammering gallery
$b \rLine^i a\rLine^j c\rLine^i d$ in $B$ (the superscripts indicate the
types of panels in the gallery). Since $ij\neq ji$, this gallery is
minimal. Therefore it is contained in some apartment $A\subseteq B$.
Let $s$ be the
panel in $A$ opposite to the $j$-panel $a\cap c$. Let $e$ be a
chamber in $S=\Res(s)$ which is not in $A$ (here we use that $B$ is thick).
Then $e$ is opposite to
$a$ and $c$. There is a unique panel $t\subseteq e$ which is opposite both to
$r$ and to the panel $q=c\cap d$. Since $b$ is not opposite $e$,
$\proj_tb=e$. Similarly $\proj_td=e$, whence $[q;t;r](b)=d$.
We claim that $\proj_t(a)=\proj_t(c)$. Assuming for the moment
that this is true,
we have $[t;r](a)=[t;q](c)$, whence $[q;t;r](a)=c$.
Applying the same construction to
$b' \rLine^i a\rLine^j c\rLine^i d$, with a second apartment $A'$ and
a panel $t'$, we obtain the projectivity $[r;t';q;t;r]$ with the
required properties.

It remains to show that $\proj_t(a)=\proj_t(c)$. Let $A''$ denote the
apartment spanned by the opposite chambers $a$ and $e$.
Let $f$ denote the chamber in $\Res(t)\cap A''$ different from $e$
and let $k$ denote the gallery distance between $a$ and $e$.
Then $a$ and $f$ have gallery distance $k-1$. Since $t$ and $p$
are not opposite, $g=\proj_pf$ has gallery distance $k-2$ from $f$.
So there is a gallery $(f,\ldots,g,a)$ of length $k-1$, which
is therefore minimal. It follows that
$(f,\ldots,g,c)$ is also a minimal gallery, and
$f$ has gallery distance $k-1$ from $c$. Since $t$ is opposite
to $r$ and $q$, we have $\proj_t(a)=f=\proj_t(c)$.
\qed
\end{Thm}
Let $a,b$ be opposite panels in $B$ and let $B(a,b)$ denote the
union of all apartments containing $a$ and $b$.
For a chamber $c\in B$, let $c_a=\proj_ac$ and $c_b=\proj_bc$.
If $A$ is an apartment containing $a$ and $b$ and $c\in A$ is a chamber,
then $c_a,c_b\in A$ and $\proj_ac_b=c_a$.
In particular,
each pair of distinct chambers $c,d\in \Res(a)$ determines
a unique apartment in $B(a,b)$ containing $c$ and $d$
and the chambers in $\Res(a)$ correspond bijectively to the
half apartments of $B(a,b)$ having $a$ and $b$ as boundary panels.
\begin{Lem}
\label{TreeSubComplex}
The subcomplex $B(a,b)\subseteq B$ is a weak building of the same type
as $B$ and every apartment of this weak building contains $a$ and $b$.
If $b'$ is another panel opposite $a$, then there is a
unique simplicial isomorphism $B(a,b)\rTo\allowbreak B(a,b')$ which fixes
$B(a,b)\cap B(a,b')$.

\proof
Let $c,d$ be chambers in $B(a,b)$. Then there exists an apartment
$A$ containing $a,b,c_a$ and $d_a$. It follows that $A$ contains
$c_b$ and $d_b$ and therefore $c$ and $d$. This shows that
$B(a,b)$ is a building and that every apartment of this building
contains $a$ and $b$.

Suppose now that $b'$ is also opposite $a$. For every apartment
$A$ containing $a,b$, there is a unique apartment $A'$
containing $a,b'$, such that $A\cap \Res(a)=A'\cap \Res(a)$.
Since $A\cap A'$ contains chambers, there is a unique isomorphism
$A\rTo A'$ fixing $A\cap A'$. For $c\in A$, the image $c'\in A'$ can
be described as follows. It is the unique chamber which has the
same $W$-valued distances from the two chambers
$a_c$ and $\proj_{b'}a_c$ as $c$ has from $a_c$ and 
$\proj_ba_c=b_c$. This description is independent of $A$ and
$A'$ and 
shows that we obtain a well-defined map on the chambers.
This map is adjacency-preserving on each apartment and hence
a building isomorphism.
\qed
\end{Lem}
\begin{Num}
\label{AutomorphismsAlongProjectivityGraph}
Suppose that $a=a_0$ is a panel and $a_0\rLine a_1\rLine\cdots\rLine a_k$ is
a path in the opposition graph. By the previous lemma we
obtain a sequence of canonical isomorphisms
$B(a_0,a_1)\cong B(a_1,a_2)\cong\ldots\cong B(a_{k-1},a_k)$
fixing the intersections of consecutive buildings.
If $a_0=a_k$, then the composition of these isomorphisms
yields an automorphism of $B(a_0,a_1)$. If $k$ is even,
then this automorphism fixes $a$
and its restriction to $\Res(a)$ coincides
with the projectivity $[a_k;\ldots;a_0]$.
\end{Num}

\section{Euclidean buildings}
\label{EuclideanBuildingsSection}
The notion of a (nondiscrete) Euclidean building is due to Tits 
\cite{TiComo}. 
Prior to their axiomatization in \cite{TiComo}, the 
nondiscrete Euclidean buildings that arise from
reductive groups over valued fields were studied in \cite{BT}.
We rely on Parreau's work \cite{Par} which contains many important
structural results for Euclidean buildings. She showed in particular
that the axioms given by Kleiner-Leeb
\cite{KL} are equivalent to Tits' original axioms plus metric
completeness.
\begin{Num}\textbf{The affine Weyl group\ }
We fix a spherical Coxeter group $(W,I)$ in its standard representation on $\RR^n$
(where $n=\#I$).
A \emph{Weyl chamber} or \emph{sector} in $\RR^n$ is the closure of a connected
component in $\RR^n-(H_1\cup H_2\cup\cdots\cup H_r)$, where
the $H_k$ are the reflection hyperplanes of $W$ corresponding to the
conjugates of the generators $i\in I$. The closure of a connected
component of $\RR^n- H_k$ is called a \emph{half space}.
A \emph{Weyl simplex} in $\RR^n$ is an intersection of Weyl chambers.
A \emph{wall} is a reflection hyperplane.
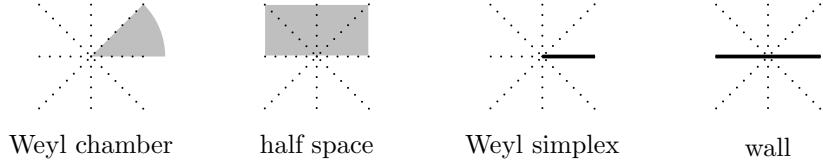
\begin{center}
\begin{pspicture}(11,2)
\psline[linewidth=0.5mm](9.3,1)(10.7,1)
\psline[linewidth=0.5mm](7,1)(7.7,1)
\psframe[linestyle=none,fillstyle=solid,fillcolor=lightgray](3.3,1)(4.7,1.7)
\pswedge[linestyle=none,fillstyle=solid,fillcolor=lightgray](1,1){1}{0}{45}
\psline[linestyle=dotted](.3,.3)(1.7,1.7)
\psline[linestyle=dotted](.3,1)(1.7,1)
\psline[linestyle=dotted](1,.3)(1,1.7)
\psline[linestyle=dotted](.3,1.7)(1.7,.3)
\psline[linestyle=dotted](3.3,.3)(4.7,1.7)
\psline[linestyle=dotted](3.3,1)(4.7,1)
\psline[linestyle=dotted](4,.3)(4,1.7)
\psline[linestyle=dotted](3.3,1.7)(4.7,.3)
\psline[linestyle=dotted](6.3,.3)(7.7,1.7)
\psline[linestyle=dotted](6.3,1)(7.7,1)
\psline[linestyle=dotted](7,.3)(7,1.7)
\psline[linestyle=dotted](6.3,1.7)(7.7,.3)
\psline[linestyle=dotted](9.3,.3)(10.7,1.7)
\psline[linestyle=dotted](9.3,1)(10.7,1)
\psline[linestyle=dotted](10,.3)(10,1.7)
\psline[linestyle=dotted](9.3,1.7)(10.7,.3)
\rput(1,-.2){\small Weyl chamber}
\rput(4,-.2){\small half space}
\rput(7,-.2){\small Weyl simplex}
\rput(10,-.2){\small wall}
\end{pspicture}\\\
\end{center}
Note that we do not require that $W$ be irreducible.
The Weyl simplices in $\RR^n$, ordered by inclusion, form a simplicial complex
which is isomorphic to the Coxeter complex $\Sigma$ of $W$.

We also fix a $W$-invariant inner product on $\RR^n$. 
The corresponding norm will be denoted by $||\cdot||$.
Up to scaling factors on the irreducible $W$-submodules of $\RR^n$,
such an inner
product is unique. The group $W$ normalizes the translation group $(\RR^n,+)$
and the semidirect product $W\RR^n$ acts isometrically on $\RR^n$. We call
this group $W\RR^n$ the \emph{affine Weyl group}.
\end{Num}
\begin{Num}\textbf{Euclidean buildings\ }
\label{EBAxioms}
Let $W$ be a spherical Coxeter group and $W\RR^n$ the corresponding
affine Weyl group.
Let $X$ be a metric space. A \emph{chart} is an isometric
embedding $\phi:\RR^n\rTo X$, and its image is called an \emph{affine apartment}.
We call two charts $\phi,\psi$
\emph{$W$-compatible} if $Y=\phi^{-1}\psi(\RR^n)$ is convex 
(in the Euclidean sense) and
if there is an element $w\in W\RR^n$ such that
$\psi\circ w|_Y=\phi|_Y$ (this condition is void if $Y=\emptyset$). 
We call a metric space $X$ together with a collection $\A$ of
charts a \emph{Euclidean building} if it has the following
properties.

\medskip\noindent
\textbf{(EB1)} For all $\phi\in\A$ and $w\in W\RR^n$, the composition
$\phi\circ w$ is in $\A$.

\medskip\noindent
\textbf{(EB2)} Any two points $x,y\in X$ are contained in some affine apartment.

\medskip\noindent
\textbf{(EB3)} The charts are $W$-compatible.

\medskip\noindent
The charts allow us to map Weyl chambers, walls and half spaces into $X$;
their images are also called Weyl chambers, walls and half spaces. The first
three axioms guarantee that these notions are coordinate independent.
We call $\A$ an \emph{atlas} or \emph{apartment system} for $X$.

\medskip\noindent
\textbf{(EB4)} If $a,b\subseteq X$ are Weyl chambers, then there is
an affine apartment $A$ such that the intersections
$A\cap a$ and $A\cap b$ contain Weyl chambers.

\medskip\noindent
\textbf{(EB5)} If $A_1,A_2,A_3$ are affine apartments which intersect pairwise
in half spaces, then $A_1\cap A_2\cap A_3\neq\emptyset$.
\end{Num}
The last axiom may be replaced by the following axiom \cite[Thm.~1.21]{Par}.

\medskip\noindent
\textbf{(EB5')} If $A\subseteq X$ is an affine apartment and $x\in A$ a point,
then there is a $1$-Lipschitz retraction $\rho:X\rTo A$ with
$\rho^{-1}(x)=\{x\}$. The point $x$ is called the \emph{center} of the retraction.

\medskip\noindent
In fact, the retractions can be chosen in such a way that they have the following
slightly stronger property.

\medskip\noindent
\textbf{(EB5$^+$)} If $A\subseteq X$ is an affine apartment and $x\in A$ a point,
then there is a $1$-Lipschitz retraction $\rho:X\rTo A$ with
$d(x,y)=d(x,\rho(y))$ for all $y\in X$. 

\medskip\noindent
The number $n$ is called the \emph{dimension} of the Euclidean building $X$.
It coincides with the topological (covering) dimension of $X$
\cite[7.1]{KramerLocal} \cite[3.3]{Lang}.
In \cite{TiComo}, \cite{KL} or \cite{Par}, the translation group of the affine
Weyl group may be some $W$-invariant subgroup of $(\RR^n,+)$.
Since we are in this paper only concerned with metric properties of Euclidean buildings,
and since the affine Weyl group can always be enlarged
to the full translation group without changing the underlying metric space and
the set of affine apartments \cite[1.2]{Par}, there is no loss in generality here.
From our viewpoint, every point $p\in X$ is a \emph{special point}
\cite[1.3.7]{BT}. We remark that the Coxeter group $W$ of a Euclidean building
need not be crystallographic \cite{HKW} \cite{BeKa}.
\begin{Num}\textbf{Example: Euclidean cones over spherical buildings\ }
\label{Cones}
Let $B$ be a spherical building and let $\EB (B)$ denote the quotient of 
$|B|\times[0,\infty)$ where $|B|\times 0$ is
identified to a point. Let $d_{|B|}$ denote the 
spherical metric on $|B|$, and put 
$d(x\times s,y\times t)^2=s^2+t^2-2st\cos(d_{|B|}(x,y))$;
see \cite[I.5.6]{BH}. With this metric,
$\EB (B)$ is the infinite Euclidean cone over $|B|$. It is not difficult to see
that $\EB (B)$ is a Euclidean building. The affine apartments in
$\EB (B)$ correspond bijectively to the apartments of $B$. These buildings
are generalizations of the trees of type (I), and we call them 
\emph{Euclidean buildings of type (I)}. One can view
$\EB (B)$ as the affine building with respect to the trivial
valuation on the spherical building $B$.
We note that this construction is functorial: every automorphism
of $B$ extends to an isometry of $\EB (B)$. In this way, every
spherical building can be viewed as a Euclidean building.
In \cite{Rou}, $\EB (B)$ is called an \emph{immeuble vectoriel}.
\end{Num}
\begin{Num}\textbf{Example: leafless trees\ }
A leafless tree is a $1$-dimensional
Euclidean building. In particular $\RR$ is a Euclidean
building. The affine Weyl group is $W\RR=\{x\mapstoo \pm x+c\mid c\in\RR\}$.
\end{Num}
\begin{Num}\textbf{Example: simplicial Euclidean buildings\ }
The geometric realization of an affine simplicial building is
a Euclidean building \cite{Bro} \cite{TiComo}.
\end{Num}
The images of the Weyl simplices in a Euclidean building
$X$ under the charts are also called Weyl
simplices. The image of the origin in $\RR^n$ is called the 
\emph{base point} or \emph{tip} of the Weyl simplex.
\begin{Num}\textbf{The vector distance\ }
\label{VectorDistance}
Let $a_0\subseteq\RR^n$ be a fixed Weyl chamber.
Given two points $p,v$ in the Euclidean building $X$, there exists a chart
$\phi:\RR^n\rTo X$ that maps $a_0$ to a $p$-based Weyl chamber
containing $v$. Let $\Theta(p,v)\in a_0\subseteq\RR^n$
denote the vector that is mapped to $v$. By (EB3), the vector
$\Theta(p,v)$ is independent of $\phi$.
We thus have a well-defined map 
\[
X\times X\rTo a_0,\qquad (p,v)\mapstoo\Theta(p,v)
\]
which we call the \emph{vector distance} \cite[1.3.1]{Par}.
We remark that the map $\Theta(p,-):X\rTo a_0$ is $1$-Lipschitz.
We also note the following. The involution $x\mapstoo -x$
on $\RR^n$ induces an involution $j:x\mapstoo w_0(-x)$ on the Weyl chamber $a_0$,
where $w_0$ is the unique longest element in the spherical Coxeter group $W$. 
We have the symmetry relation
\[
 \Theta(x,y)=j(\Theta(y,x)),
\]
for all $x,y\in X$.
\end{Num}
The affine apartment $X=\RR^n$ with $\A=W\RR^n$ is an example of a Euclidean building and
we have the following small fact which will be needed later. A related result
is proved in \cite[2.15]{Par}.
\begin{Lem}
\label{LocalCompatibilityCriterion}
Let $W$ be a spherical Coxeter group, acting in the standard representation on
$\RR^n$. Let $\Theta$ denote the corresponding vector distance on $\RR^n$.
Suppose that $K\subseteq \RR^n$ is a nonempty convex set, that $g:K\rTo K$ is a bijection,
and that 
\[
\Theta(p,v)=\Theta(g(p),g(v))
\]
holds for all $p,v\in K$.
Then there exists an element $w$ in the affine Weyl group  $W\RR^n$ whose restriction to $K$ 
coincides with $g$.

\proof
First we note that $g$ is an isometry of $K$, because $||\Theta(p,v)||=d(p,v)$.
We choose the pair $p,v\in K$ in such a way that the smallest Weyl simplex $b\subseteq a_0$
that contains $\Theta(p,v)$ has maximal dimension. Applying an element of
the affine Weyl group $W\RR^n$ to $K$, we may assume that $p=0$ and that $v\in b\subseteq a_0$.
We choose $w\in W\RR^n$ in such a way that $w(0)=g(0)$ and $w(v)=g(v)$. We claim
that $g(x)=w(x)$ for all $x\in K$.

In general, if $\Theta(0,u)=\Theta(0,u')$, then $u$ and $u'$ lie in the same
finite $W$-orbit in $\RR^n$. Since $b$ is maximal,
$v$ is an interior point of the Weyl simplex $b$ and
there is an $\eps>0$ such that $g(u)\in w(b)$ holds for all $u\in b\cap K$
with $d(u,v)\leq\eps$. This in turn implies that $g(u)=w(u)$ for all $u\in b\cap K$
with $d(u,v)\leq\eps$. Since both $g$ and $w$ are isometries, we conclude that
$g$ and $w$ agree on the convex set $b\cap K$.
Let now $u\in K$ be arbitrarily. If $t>0$ is small we have necessarily
$z=ut+v(1-t)\in b$ by the maximality of $b$. Since both $g$ and $w$ are isometries
and $g(v)=w(v)$ and $g(z)=w(z)$, we have $g(u)=w(u)$.
\qed
\end{Lem}
The Weyl simplices lead to two spherical buildings.
One of them captures the asymptotic geometry of $X$, while the other
is an `infinitesimal' version of the Euclidean building, similar to
the tangent space of a Riemannian manifold. The first one,
the spherical building at infinity, will be
considered now and the second in \ref{ResidueAtPoint}.
\begin{Num}\textbf{The spherical building at infinity\ }
We call two Weyl simplices $a,a'\subseteq X$ \emph{Hausdorff equivalent} if they
have finite Hausdorff distance. The equivalence class of
$a$ is denoted $\partial a$. The preorder $\subseteq_{Hd}$
defined in \ref{FinHd} induces a partial order on these
equivalence classes. Let
$\partial_{\A} X$ denote the set of all equivalence classes
of Weyl simplices, partially ordered by
domination $\subseteq_{Hd}$. For every affine apartment $A$,
the poset $\partial A$ consisting of the
Weyl simplices contained in $A$ may be viewed as a sub-poset of
$\partial_{\A} X$. 
\end{Num}
\begin{Prop}
The poset $\partial_{\A}X$ is a spherical building. 
The map $A\mapstoo\partial A$
is a one-to-one correspondence between the affine apartments
in $X$ and the apartments of the spherical building
$\partial_{\A}X$.

\proof
See Parreau \cite[1.5]{Par}.
\qed
\end{Prop}
We have also the following fact.
\begin{Lem}
\label{ConvexityLemma}
Suppose that $a$ and $a'$ are Weyl chambers with tips $x,x'$.
If $\partial a=\partial a'$, then 
\[
Hd(a,a')\leq d(x,x').
\]

\proof
Let $\psi:a\rTo a'$  denote the unique $\Theta$-invariant isometry. Given
$z\in a$, there is a unique 
geodesic $\gamma:[0,\infty)\rTo a$ with $\gamma(0)=x$ and $\gamma(s)=z$, for some $s\geq0$.
The geodesic $\gamma'=\psi\circ\gamma$ is at finite distance from $\gamma$,
because $\partial a=\partial a'$.
Since the metric of a CAT$(0)$ space is convex, 
we have $d(\gamma(s),\gamma'(s))\leq d(x,x')$ \cite[II.2.2]{BH}.
\qed
\end{Lem}
\begin{Num}\textbf{The maximal atlas\ }
\label{MaximalAtlas}
The spherical building at infinity depends very much on the
chosen set of charts $\A$.
Similarly to a differentiable manifold,
a Euclidean building admits a unique maximal atlas $\Amax$
\cite[2.6]{Par}.
The set $\Amax$ is called the \emph{complete apartment system}.
We denote the spherical building at infinity which corresponds
to the complete apartment system by $\partialmax X$.
Its metric realization coincides with the Tits boundary of
$X$ \cite[Cor.~2.19]{Par}.
\end{Num}
\begin{Num}\textbf{Product decompositions and reductions\ }
\label{ThickReductionEuclideanBuilding}
If the Coxeter group $W$ is reducible, then $X$ decomposes
as a metric product $X=X_1\times X_2$ of Euclidean buildings,
with 
$\partial_{\A_1}X_1*\partial_{\A_2}X_2
=\partial_{\A}X$ \cite[2.1]{Par}.
If $\SS^k*B_0=\partial_{\A}X$ is
the thick reduction of $\partial_{\A}X$ and if
$W_0$ is the Coxeter group of $B_0$, then there is a Euclidean
building $X_0$ with affine Weyl group $W_0\RR^m$ and atlas
$\A_0$ and an isometry
$X\cong\RR^{n-m}\times X_0$, with $\partial_{\A_0}X_0=B_0$
and $k=n-m-1$ \cite[4.9]{KL}.
\end{Num}
We record some more results about Euclidean buildings which
will be needed later.
\begin{Prop}\label{MaximalFlatsProp}
Let $X$ be an $n$-dimensional Euclidean building. If
$A\subseteq X$ is isometric to $\RR^n$, then $A$ is an
apartment in the complete apartment system of $X$.

\proof
See \cite[Prop.~2.25]{Par}.
\qed
\end{Prop}
We remark that the previous proposition is in the approach by Kleiner and Leeb
essentially an axiom.
\begin{Num}\textbf{The residue at a point\ }
\label{ResidueAtPoint}
Let $p$ be a point in $X$. Two $p$-based Weyl simplices $a,b$
have \emph{equal germs near $p$} if 
\[
B_r(p)\cap a=B_r(p)\cap b
\]
holds for some $r>0$. The corresponding equivalence classes
form a spherical building $X_p$ which we
call the \emph{residue} of $X$ at $p$; see  \cite[1.6]{Par}.
The Coxeter group of $X_p$ is $W$, hence $X_p$ has the same type
as the spherical building at infinity $\partial_\A X$.
We call the point $p$ \emph{thick} if the residue $X_p$ is thick.
If $a\subseteq X$ is a Weyl simplex, then $\partial a$ has
a unique representative which is a $p$-based Weyl simplex
\cite[Cor.~1.9]{Par}. The apartments of $X_p$ arise from
the affine apartments containing $p$ (but this correspondence is
in general not $1$-to-$1$). We thus obtain 
a type-preserving surjective simplicial map
\[
\partial_{\A}X\rTo X_p,
\]
and also a surjective canonical
$1$-Lipschitz map $\EB (\partial_{\A}X)\rTo X$ (which depends on $p$).
\end{Num}
\begin{Num}
\label{AngelRigidity}
Given a $p$-based Weyl chamber $a$, let $\zeta_a$ denote the spherical
barycenter of the Weyl chamber $\partial a$. We call the geodesic ray
\[
[p,\zeta_a)\subseteq a
\]
the \emph{central ray} of $a$.
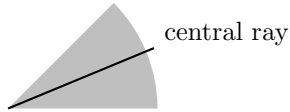
\begin{center}
\begin{pspicture}(0,1)(4,2.5)
\pswedge[linestyle=none,fillstyle=solid,fillcolor=lightgray](1,1){2}{0}{45}
\psline[linestyle=solid](1,1)(2.940,1.804)
\rput(3.9,2){\small central ray}
\end{pspicture}
\end{center}
The possible angles between all centrals rays
starting at $x$ form a finite subset of $[0,\pi]$, since the 
spherical Coxeter complex of $W$ is finite. If $[p,\zeta_b)$ and
$[p,\zeta_c)$ are the central rays of two $p$-based Weyl chambers
$b,c$, then the Alexandrov angle between these rays is $0$ if and
only if $b$ and $c$ have the same germ near $p$.
\end{Num}
If two sufficiently long geodesic segments in a tree have
small Hausdorff distance, then they intersect in a long
geodesic segment. The next proposition is
the higher dimensional analog of this fact.
\begin{Prop}
\label{CloseApartmentsIntersect}
Let $X$ be a Euclidean building and let $s>0$ be a positive real number.
There exists a positive real number $c_s>0$ 
such that the following holds for all $t>0$.
If $A,A'\subseteq X$ are affine apartments and $y$ is a point in $X$
such that $Hd_{ B_{t+c_s}(y)}(A,A')\leq s$, then 
\[
B_t(y)\cap A=B_t(y)\cap A'.
\]

\proof
Let $a_0\subseteq\RR^n$ denote the standard Weyl chamber.
There exists a number $c_s>0$ such that for the point 
$z_a\in[0,\zeta_a)$ with $d(0,z_a)=c_s$, the closed
ball $\bar B_{2s}(z_a)\cap a$ contains only interior points of $a$. 

Suppose now that $t>0$, that $A,A'\subseteq X$ are as in the claim 
of the proposition and that $p\in B_t(y)\cap A$. We claim that
$p\in B_t(y)\cap A'$. 
We choose a $p$-based Weyl chamber $b\subseteq A$
and $z_b\in[p,\zeta_p)$ with $d(p,z_b)=c_s$.
The condition on the Hausdorff distance ensures that we can find
a point $v'\in A'$ with $d(z_b,v')\leq 2s$.
From the comparison triangle in Euclidean space we see that
the angle between the geodesic segments $[p,v']$ and $[p,z_b]$ is
so small that $[p,v']$ intersects the Weyl chamber $b$ in an interior point.
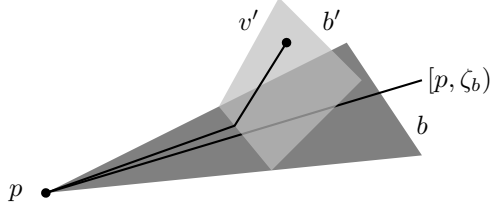
\begin{center}
\begin{pspicture}(4,3)
\pspolygon[linestyle=none,fillstyle=solid,fillcolor=gray](0,0)(5,0.5)(4,2)
\psline[linestyle=solid](0,0)(5,1.5)
\pspolygon[linestyle=none,fillstyle=solid,fillcolor=lightgray,opacity=0.7](3,0.3)(4.2,1.5)(3.1,2.6)(2.3,1.15)
\psline[linestyle=solid]{*-*}(0,0)(2.51,0.895)(3.2,2)
\rput(-.4,0){\small$p$}
\rput(2.7,2.3){\small$v'$}
\rput(5.5,1.5){\small$[p,\zeta_b)$}
\rput(3.8,2.3){\small$b'$}
\rput(5.0,.9){\small$b$}
\end{pspicture}
\end{center}
We now choose a geodesic ray $[v',\eta)\subseteq A'$ which extends
$[p,v']$. Thus $[p,v']\cup[v',\eta)$ is a geodesic ray. 
This ray is contained in a unique $p$-based Weyl chamber $b'$ whose germ
near $p$ is contained in $A$.

We now consider the $p$-based Weyl chamber $c\subseteq A$ opposite
$b$, and we choose $z_c\in[p,\zeta_c)$ similarly as before, with
$d(p,\zeta_c)=c_s$. Then we find a $p$-based
Weyl chamber $c'$ whose germ near $p$  is opposite to $b'$. 
These two Weyl chambers are contained in a unique affine apartment
\cite[Prop.~1.12]{Par}, which is therefore
$A'$. It follows that $p$ is contained in $A'$. Thus
$B_r(y)\cap A\subseteq B_r(y)\cap A'$. The other inclusion follows by symmetry.
\qed
\end{Prop}
We call an affine wall $M\subseteq X$
\emph{thick} if $M$ is the intersection of three affine apartments.
We have the following local thickness criterion for walls.
\begin{Lem}
\label{ThickWallLemma}
Let $A\subseteq X$ be an affine apartment, let $p\in A$ and
let $M\subseteq A$ be a wall containing $p$. Let $r\in X_p$
be a panel in the wall determined by $M$ in $X_p$.
If $r$ is contained in three distinct chambers of $X_p$, then $M$ is thick.

\proof
Let $A_p$ denote the apartment induced by $A$ in $X_p$, let
$s\in A_p$ be the panel opposite $r$ and let
$a,b\in A_p$ be the two chambers containing $s$.
We represent $a,b,r,s$ by $p$-based Weyl simplices in 
$A$. By assumption, there is a chamber containing $r$ which is
opposite both
to $a$ and to $b$. By \ref{ThicknessPullsBack} we can find a $p$-based
Weyl chamber $c$ representing this chamber, such that
$\partial c$ is adjacent to $\partial A$. 
Since the panel $\partial c\cap\partial A$ has a unique
$p$-based representative, this representative is contained
in $M$. It follows that $c\cap A\subseteq M$.
Let $A'$ denote the affine apartment spanned by
$c$ and $a$; see \cite[1.12]{Par}.
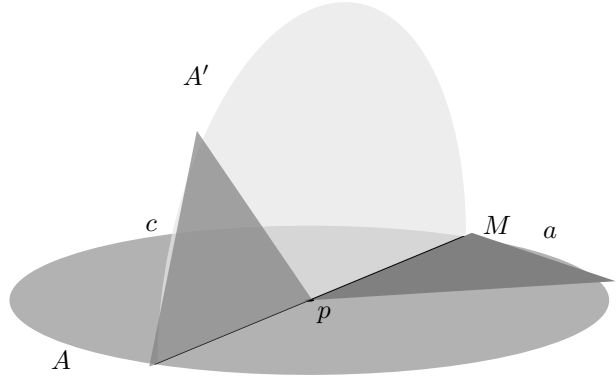
\begin{center}
\begin{pspicture}(-2,-1.0)(4,4)
\psset{viewpoint=-5 -3 1.5}
\ThreeDput[normal=0 0 1]{
\pscircle[linestyle=none,fillstyle=solid,fillcolor=gray7](0,0){4}
\pspolygon[linestyle=none,fillstyle=solid,fillcolor=gray](0,0)(4.2,0)(2.969,-2.969)
\psline[linestyle=solid](-4,0)(4,0)
\psline[linestyle=solid]{*-*}(0,0)(0,0)
}
\ThreeDput[normal=0 -1 0](0,0,0){
\psclip{\psframe[linestyle=none](-5,0)(5,5)}
\pscircle[linestyle=none,fillstyle=solid,fillcolor=gray9,opacity=.7](0,0){4}
\endpsclip
\pspolygon[linestyle=none,fillstyle=solid,fillcolor=gray,opacity=0.7](0,0)(-4.2,0)(-2.969,2.969)
}
\rput(.2,-.2){\small$p$}
\rput(2.5,1.0){\small$M$}
\rput(3.2,0.9){\small$a$}
\rput(-2.1,1){\small$c$}
\rput(-3.3,-.8){\small$A$}
\rput(-1.5,3.0){\small$A'$}
\end{pspicture}
\end{center}
Then $A\cap A'$ is a half space whose boundary is
$M$ and therefore $M$ is thick.
\qed
\end{Lem}

\begin{Lem}
\label{ThickPointLem}
Let $p$ be a point of an affine apartment $A\subseteq X$. Then  
$p$ is thick if and only if  every wall of $A$ containing $p$ is
thick.

\proof
If $p$ is thick, then every wall through $p$ is thick by
\ref{ThickWallLemma}.
Conversely, if every wall of $A$ containing
$p$ is thick, then all panels in $A_p$ are contained in at least
$3$ chambers. Thus, $p$ is thick by \ref{ThickApartmentsLemma}.
\qed
\end{Lem}
\begin{Lem}
\label{ThickWallReflection}
Let $A$ be an affine apartment and assume that $L,M\subseteq A$
are non parallel thick walls. Then the reflection of $L$
along $M$ in $A$ is also thick.

\proof
Let $A_\pm \subseteq A$ denote the half spaces bounded by $M$ and let
$A_0\subseteq X$ be a third
half space with $A\cap A_0=M$. There is a unique
affine wall $L'$ in the affine apartment $A'=A_+\cup A_0$ 
extending $L\cap A_+$. By \ref{ThickWallLemma},
$L'$ is thick. Now let $L''$ be the wall in
$A_0\cup A_-$ which extends $L'\cap A_0$.
Again by \ref{ThickWallLemma}, $L''$ is thick. Finally,
let $L'''$ be the wall in $A$ which extends $L''\cap A_-$.
A third application of \ref{ThickWallLemma} yields
that $L'''$ is thick. 
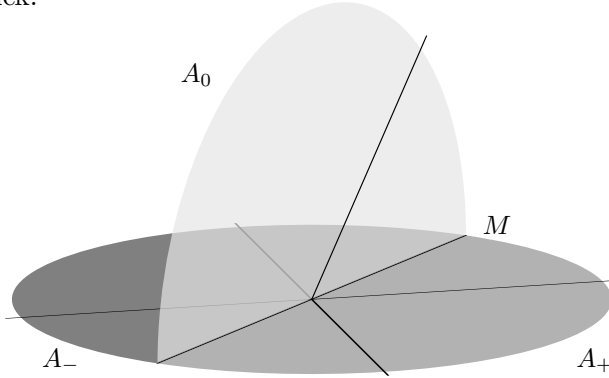
\begin{center}
\begin{pspicture}(-2,-1.0)(4,3.5)
\psset{viewpoint=-5 -3 1.5}
\ThreeDput[normal=0 0 1]{
\psclip{\psframe[linestyle=none](-5,-5)(5,0)}
\pscircle[linestyle=none,fillstyle=solid,fillcolor=gray7](0,0){4}
\endpsclip
\psclip{\psframe[linestyle=none](-5,0)(5,5)}
\pscircle[linestyle=none,fillstyle=solid,fillcolor=gray5](0,0){4}
\endpsclip
\psline[linestyle=solid](-2.969,-2.969)(0,0)(-2.969,2.969)
\psline[linestyle=solid](2.969,-2.969)(0,0)(2.969,2.969)
\psline[linestyle=solid](-4,0)(4,0)
}
\ThreeDput[normal=0 -1 0](0,0,0){
\psclip{\psframe[linestyle=none](-5,0)(5,5)}
\pscircle[linestyle=none,fillstyle=solid,fillcolor=gray9,opacity=.7](0,0){4}
\endpsclip
\psline[linestyle=solid](2.969,2.969)(0,0)
}
\rput(2.5,1.0){\small$M$}
\rput(-3.3,-.8){\small$A_-$}
\rput(3.8,-.8){\small$A_+$}
\rput(-1.5,3.0){\small$A_0$}
\end{pspicture}
\end{center}
But $L'''$ is precisely the reflection of $L$ along $M$.
\qed
\end{Lem}
\begin{Lem}
\label{ThickPointsExist}
The Euclidean building $X$ contains a thick point if and only
if $\partial_\A X$ is thick.

\proof
By \ref{ThicknessPullsBack} thickness of $X_p$ implies thickness of
$\partial_\A X$.
Now suppose that $\partial_\A X$ is thick, let
$A$ be an affine apartment, let $c$ be a Weyl chamber
of $A$ and let $M_1,...,M_n$ be the walls of $A$
bounding $c$. Since $\partial_\A X$ is thick,
we can choose thick walls $M_1',\ldots,M_n'$
in $A$ such that $M_i'$ is parallel to $M_i$ for each $i$.
The intersection of the thick 
walls $M_1',\ldots,M_n'$ contains a
point $p$. By \ref{ThickPointLem} and \ref{ThickWallReflection}, $p$ is thick.
\qed
\end{Lem}
The following trichotomy is analogous to the case of trees in \ref{TreeStructureCor}.
In higher dimension, however, we need no assumptions on group actions.
\begin{Prop}
\label{Trichotomy}
Suppose that $X$ is a Euclidean building of dimension $n\geq 2$
and that $\partial_\A X$ is irreducible and thick. Let
$th(X)\subseteq X$ denote the set of thick points.
There are the following three possibilities.

(I) There is a unique thick point which is contained in every affine apartment
of $X$.

(II) The set of thick points is a closed, discrete and cobounded
subset in $X$ and in  every apartment of $X$.

(III) The set of thick points is dense in $X$ and in every apartment of $X$.

\proof
We noted in \ref{ThickPointsExist} that $th(X)\neq\emptyset$.
For an affine apartment $A\subseteq X$ we denote by $R(A)\subseteq \Isom(A)$ the
group generated by reflections along the thick walls in $A$.
If $p$ is a thick point, then the $R(A)$-stabilizer of $p$ is $R(A)_p\cong W$.
We start with two observations.

(i)
Suppose that $p$ is a thick point in an affine apartment $A$,
and that $M\subseteq A$ is a thick wall not containing $p$.
Such an apartment $A$ exists if $X$ contains at least two
thick points.
Since the Weyl group $W$ is irreducible and $\dim(A)\geq 2$,
the point $p$ is the
intersection of thick walls in $A$ which are not parallel
to $M$. It follows from \ref{ThickWallReflection} 
that the reflection of $p$ along $M$ in $A$
is again a thick point, so the $R(A)$-orbit of $p$
consists of thick points.

(ii) For any two affine apartments $A,A'$, 
there is a sequence of affine apartments
$A=A_0,\ldots,A_k=A'$ such that $A_j\cap A_{j+1}$ is
a half space. (Using galleries, this is easily seen to be true
for apartments in the spherical building at infinity.)
It is clear that the sequence of walls determined by these
half spaces is thick.

Suppose first that $th(X)=\{p\}$ consists of a single point and
that $A$ is an affine apartment containing $p$. Then (i)
shows that the thick walls in $A$ are precisely
the ones passing through $p$. If $A'$ is any other affine apartment,
then (ii) and a simple induction on $k$ show that all apartments
$A_0,\ldots,A_k$ in the sequence connecting $A$ and $A'$ 
contain $p$. This is case (I).

Suppose now that $th(X)$ contains at least two points $p,q$
and that $A$ is an apartment containing $p$ and $q$.
Then the $R(A)$-orbit of $q$ consists of thick points and
is cobounded in $A$, because the Weyl group $W$
is irreducible, and because $R(A)$ contains translations.
Let $T\subseteq R(A)$ denote the translational part
(i.e. the kernel of the action of $R(A)$ on $\partial A$).
Any closed subgroup of the vector group $(\RR^n,+)$ is
a product of a free abelian group of finite rank and a vector group.
If $T$ is discrete, then $R(A)$ is an affine reflection
group. If $T$ is not discrete, then
the closure of $T$ in $\Isom(A)$ consists of all translations
in $A$ (because $W$ acts irreducibly on the set of all translations),
so $th(X)\cap A$ is dense.

If $A'$ is any other
affine apartment having a half space in common with $A$,
then the thick walls in $A$ propagate to thick walls
in $A'$. The isometry $A\rTo A'$ also preserves thick points.
We obtain a canonical isomorphism
$R(A)\cong R(A')$. From (ii) we conclude that the
isometry groups $R(A)$ and $R(A'')$ are isomorphic for
any affine apartment $A''\subseteq X$, and that this
isomorphism maps $th(X)\cap A$ onto $th(X)\cap A''$.
In the discrete case, we have (II), whereas the
nondiscrete case is (III).
\qed
\end{Prop}
\begin{Cor}
Let $X$ be as in \ref{Trichotomy} (I). Then $X\cong \EB (\partial_\A X)$
is a Euclidean cone.

\proof Let $p\in X$ be the unique thick point and let
$A$ be an affine apartment. Since we are in case (I),
the point $p$ is contained in $A$. The canonical surjective
$1$-Lipschitz map $\EB (\partial_{\A}X)\rTo X$ sends therefore
$\EB (\partial A)$ isometrically onto $A$. Since any two points
in $\EB (\partial_{\A}X)$ are in some apartment,
$\EB (\partial_{\A}X)\rTo X$ is an isometry.
\qed
\end{Cor}
The following observation will not be used, but it illustrates
how simplicial affine buildings fit into the picture.
In the remaining case (III), $X$ is a nondiscrete Euclidean building.
\begin{Cor}
Let $X$ be as in \ref{Trichotomy} (II). Then $X$ is the metric
realization of an affine simplicial building.

\proof
The groups $R(A)$ are affine Coxeter groups; see \cite[Ch.VI]{Bro}.
In this way, every apartment has a canonical simplicial structure
as a Coxeter complex. The axioms of an affine simplicial building
follow from \ref{EBAxioms}.
\qed
\end{Cor}
In general, a Euclidean building is not determined by its
spherical building at infinity. For example,
\[
\partial_\A \EB (\partial_\A X)=\partial_\A X.
\]
Some additional data are needed,
which are encoded in the panel trees. The following material 
is due to Tits \cite{TiComo}. A proof of \ref{TitsEcologicalTheorem} with
all details filled in can be found in \cite{Hitz}.
\begin{Num}\textbf{Wall trees and panel trees\ }
\label{WallTrees}
Let $(X,\A)$ be a Euclidean building and let $a,b$ be a
pair of opposite panels in the spherical building at infinity.
Let $X(a,b)$ denote the union over all
affine apartments in $\mathcal A$ whose boundary contains $a$ and $b$,
\[
 X(a,b)=\bigcup\{A\mid A\subseteq X \text{ affine apartment in }\A\text{ and }a,b\in\partial A\}.
\]
Then $X(a,b)$ is a Euclidean building, which factors metrically
as 
\[
X(a,b)=T\times\RR^{n-1},
\]
where $T$ is a leafless tree;
see \cite[4.8.1]{KL}, \cite[3.9]{Rou2} \cite{TiComo}.
We call this tree $T$ the \emph{wall tree} associated to $(a,b)$,
because it depends only on the unique wall of 
$\partial_\A X$ containing $a$ and $b$.
In the notation of \ref{TreeSubComplex}, the spherical
building at infinity of $X(a,b)$ is 
\[
\partial(X(a,b))=(\partial_\A X)(a,b).
\]
If $X$ is metrically complete, then $T$ and $X(a,b)$ are also metrically
complete \cite{StruyveComplete}.
\end{Num}
\begin{Lem}
Let $X$ be a Euclidean building and let $a,b,b'$ be panels in
$\partial_\A X$. Suppose that $b$ and $b'$ are opposite $a$.
Then there is a unique isometry $X(a,b)\rTo X(a,b')$ which fixes
$X(a,b)\cap X(a,b')$ pointwise.

\proof
Let $A$ be an affine apartment whose boundary contains $a,b$,
and let $A'$ be the corresponding affine apartment whose
boundary contains $\Res(a)\cap\partial A$ and $b'$; see \ref{TreeSubComplex}.
Since $A\cap A'$ contains Weyl chambers, there is a unique
isometry $A\rTo A'$ fixing $A\cap A'$. This proves the
uniqueness of the isometry. For the existence, we show that
these isometries $A\rTo A'$ of the individual
affine apartments fit together.

Let $A_1,A_2\subseteq X(a,b)$ be two affine apartments containing
a point $x\in A_1\cap A_2$. Correspondingly, we have chambers
$c_i,d_i\in \Res(a)$ such that $c_i,\proj_b d_i$ spans $\partial A_i$;
see \ref{TreeSubComplex}. If two of these four chambers coincide, say
$c_1=c_2$, then $A_1$ and $A_2$ have the $x$-based Weyl chamber $\tilde c$
representing $c_1$ in common.
Then $A_1,A_2,A_1'$ and $A_2'$ have a sub-Weyl chamber of $\tilde c$
in common, and therefore the two isometries $A_1\rTo A_1'$
and $A_2\rTo A_2'$ map $x$ to the same point $x'$.

If $c_1,c_2,d_1,d_2$ are pairwise different, then $x$ is also
in the affine apartment determined by $c_1,d_2$ or $c_1,c_2$, because
$X(a,b)$ is a product of a tree and Euclidean space.
The previous argument, applied twice, shows that the various
apartment isometries coincide on $x$. Thus we have a well-defined
bijection which is apartment-wise an isometry. Since any two points
are in some affine apartment, the map is an isometry.
\qed
\end{Lem}
\begin{Num}
\label{WallTreesAndProjectivities}
If $b'$ varies over the
panels opposite $a$, we obtain a family of trees which are
pairwise canonically isomorphic. This canonical
isomorphism type of a tree is the \emph{panel tree} $T_a$
associated to $a$.

If $a_0\rLine a_1\rLine\cdots\rLine a_k$ is a path
in the opposition graph, then the isomorphisms
$(\partial_\A X)(a_0,a_1)\cong(\partial_\A X)(a_1,a_2)\cong\cdots$
from \ref{TreeSubComplex} are accompanied by isometries
$X(a_0,a_1)\cong X(a_1,a_2)\cong\cdots$. In particular, the
group of projectivities $\Pi_{\partial_A X}(a)$
acts on the wall tree and on the panel tree $T_a$.
If we restrict to even projectivities, the action on the
Euclidean factor is trivial.
If the building is thick and the type of the panel $a$ is
not isolated, then the action of $\Pi^+_{\partial_\A X_1}(a)$
on the ends of $T_a$ is
$2$-transitive by \ref{KnarrThm}

The branch points in the panel trees correspond to the
thick walls in $X$. Assuming that $\partial_\A X$
is irreducible, we have the following
consequence of \ref{Trichotomy}.
If one panel tree $T_a$ is of type (I), then every
panel tree is of type (I) and
$X$ is a Euclidean cone over $|\partial_\A X|$ as
in \ref{Cones}.
If one panel tree is of type (II), then every
panel tree is of type (II) and $X$ is simplicial.
The remaining possibility is that every panel tree
is of type (III).
\end{Num}
\begin{Num}\textbf{Ecological boundary isomorphisms\ }
Let $X_1,X_2$ be irreducible Euclidean buildings. Assume that
$\partial_{\A_1}X_1$ and $\partial_{\A_2}X_2$ are
thick and that 
\[
\phi:\partial_{\A_1}X_1\rTo\partial_{\A_2}X_2
\]
is an isomorphism.
Assume moreover that for every panel $a\in \partial_{\A_1}X_1$,
there is a tree isometry 
\[
\phi_a:T_a\rTo T_{\phi(a)}.
\]
If for each panel $a$, the map 
$\partial\phi_a:\partial T_a\rTo\partial T_{\phi(a)}$
coincides with the restriction of $\phi$ to $\Res(a)$, then
$\phi$ is called \emph{tree-preserving} or \emph{ecological}.
The following is \cite[Thm.~2]{TiComo}, see also \cite[Ch.~7]{Hitz}.
\end{Num}
\begin{Thm}{\rm\bf(Tits)}
\label{TitsEcologicalTheorem}
Let $(X_1,\A_1)$ and $(X_2,\A_2)$ be Euclidean buildings.
Assume that $\partial_{\A_1} X_1$ is thick and irreducible.
If $\phi:\partial_{\A_1}X_1\rTo^\cong\partial_{\A_2}X_2$ is an
ecological isomorphism,
then $\phi$ extends to an isomorphism
$\Phi:(X_1,\A_1)\rTo^\cong(X_2,\A_2)$.
\qed
\end{Thm}
The irreducibility is actually not important for the proof,
but the result is stated in this way in \cite{TiComo}.
For our purposes, the irreducible version suffices.
We remark that the following interesting problem is open.
\begin{Num}\textbf{Conjecture\ }
{\em Let $X_1,X_2$ be Euclidean buildings of type (II)
and dimension $n\geq 2$. Assume that $\partial_{cpl}X_1$ is irreducible.
Then every isomorphism 
\[
 \partial_{cpl}X_1\rTo\partial_{cpl}X_2
\]
extends uniquely to an isometry $X_1\rTo X_2$.}
\end{Num}
This is known to be true if $\partial_{cpl}X_1$ is Moufang
\cite[27.6]{WeAff}. The proof is algebraic and uses
the fact that a field admits at most one discrete complete
valuation. The conjecture is false for Euclidean buildings of
type (III), since for example $\CC$ admits infinitely many nonisomorphic
nondiscrete complete valuations, corresponding to nonisomorphic
Euclidean buildings of type $\widetilde{\mathsf{A}}_m$. The boundary is
always the spherical building of $\mathrm{SL}_{m+1}\CC$.

\section{Coarse equivalences of Euclidean buildings}

In this final section we prove, among other things, Theorem II and III
from the introduction. An important technical tool is the 'higher dimensional
Morse Lemma' \ref{HigherDimensionalMorseLemma}.
The main step in this result, in turn, is the
following topological rigidity theorem \cite[6.4.2]{KL} \cite[Sec.~4]{KramerLocal}. 
\begin{Thm}[Topological rigidity of affine apartments]\label{TopologicalRigidityOfApartments}
Let 
\[
f:X\rTo Y
\]
be a homeomorphism of metrically complete Euclidean
buildings. If $A\subseteq X$ is an affine apartment, then
$f(A)$ is an affine apartment in the complete apartment system
of $Y$.
\skop
We give an outline of the sheaf-theoretic proof from \cite[Sec.~6]{KramerLocal}.
On any Hausdorff space $X$, there are the following two  presheaves.
The \emph{local homology presheaf} assigns to every open subset $U\subseteq X$
the singular homology group $H_*(X,X\setminus U)$.
The stalk of the corresponding sheaf $\mathcalligra H\,_*(X)$
at a point $p\in X$ is the local homology group 
\[
\mathcalligra H\,_*(X)_p\cong H_*(X,X\setminus\{p\}).
\]
The second presheaf assigns to $U$ the poset 
$Cls(U)=\{A\cap U\mid A\subseteq X\text{ is closed }\}$ of all relatively 
closed subsets of $U$. 
Its stalk $\mathcalligra{Cls\,}_p$ consists of germs of closed sets at $p$.
There is a natural
transformation $supp$ between these (pre)sheaves which assigns to every relative
cycle $c\in H_*(X,X\setminus U)$ its support 
\[
supp(c)=\{p\in U\mid
c \text{ has nontrivial image in }H_*(X,X\setminus\{p\})\}
\]
\cite[6.2]{KL} \cite[Sec.~4]{KramerLocal}.
If $X$ is a Euclidean building, then the residue $X_p$ is by definition
a sub-poset of $\mathcalligra{Cls\,}_p$. The main step in the proof 
of \ref{TopologicalRigidityOfApartments} is to show
that this sub-poset can be described by means of the transformation $supp$.
The reason is as follows. For any CAT$(0)$ space $X$,
the natural map  
\[
X\setminus\{p\}\rTo\Sigma_pX
\]
which assigns to a
point $q\neq p$ the direction of the geodesic segment $[q,p]$ in the
(completed) \emph{space of directions} $\Sigma_pX$ is a homotopy equivalence
\cite[Thm.~A]{KramerLocal}. In particular, there is a natural isomorphism
in homology
\[
H_{*+1}(X,X\setminus\{p\})\rTo^\cong \tilde H_*(\Sigma_pX).
\]
If $X$ is a Euclidean building, then the space of directions
$\Sigma_pX$ is the geometric realization of the residue $X_p$.
The reduced homology (of the geometric realization) of the spherical building $X_p$
is then a free $\ZZ$-module spanned by all apartments containing a given
chamber, see \cite[Sec.~5]{KramerLocal}, the \emph{Steinberg module}
(this is essentially the contents of the Solomon-Tits Theorem).
In particular, the dimension $n=\dim(X)$ can be read off from
the local homology groups (alternatively, the number $n$ coincides with
the covering dimension of $X$ \cite[7.1]{KramerLocal} \cite[3.3]{Lang} 
and hence is a topological invariant).
This allows us to describe the image of the transformation 
\[
 \mathcalligra H\,_*(X)_p\rTo^{\ supp\ }\mathcalligra{Cls,}_p.
 \]
If the residue $X_p$ is thick, then the simplices
of $X_p$ are precisely the indecomposable elements in the distributive
poset lattice generated by $supp(\mathcalligra H\,_*(X)_p)$, see 
\cite[Cor.~6.6]{KramerLocal} (we call an element indecomposable if it is
not a union of finitely many strictly smaller elements).
This is clearly a topological invariant of the space $X$.
In general, we have by \ref{ThickRedSph} a thick reduction $X_p\cong\SS^k*|B_0|$,
with a thick spherical building $B_0$.
In this case the indecomposable elements in the image of $supp$
correspond to the simplices of the thick part
$B_0$, see \cite[5.6]{KramerLocal}. In particular, we can read off
the number $k$ from the difference between the dimension of the simplicial
complex $B_0$ and the degree $j$ for which $\tilde H_j(X\setminus\{p\})\neq 0$.

It follows that the homeomorphism maps a small open ball 
$B_r(p)\cap A$ in the affine apartment $A\subseteq X$ to a small open set
in some affine apartment in $Y$. Thus $f(A)\subseteq Y$ is a complete
simply connected flat Riemannian manifold, and therefore
isometric to $\RR^n$. By \ref{MaximalFlatsProp},
$f(A)$ is an affine apartment in the complete apartment system of $Y$.
See \cite[Sec.~6]{KramerLocal} for details and a more general result.
If $\partial_\A X$ is thick and irreducible and if $n\geq 2$, then the
completeness assumption can be dropped.
\qed
\end{Thm}
We also have the following corollary of the proof.
\begin{Cor}
\label{HomeosPreserveEuclideanFactors}
Assume that $X$ and $Y$ are metrically complete Euclidean
buildings, and that $\partialmax X$
and $\partialmax Y$ are thick.
If $f:X\times\RR^{m_1}\rTo Y\times\RR^{m_2}$
is a homeomorphism, then $m_1=m_2$.

\proof
Let $p\in X$ be a point. Then the residue at $p\times v\in X\times\RR^{m_1}$
decomposes as 
\[
(X\times\RR^{m_1})_{p\times v}=|B_0|*\SS^k,
\]
where
$B_0$ is a thick spherical building and $k\geq m_1$. If $p$ is thick, then
$k=m_1$. The same applies to $Y$, hence $m_1=m_2$.
\qed
\end{Cor}
Theorem \ref{TopologicalRigidityOfApartments} is the `hard part'
of the proof of Theorem \ref{HigherDimensionalMorseLemma} below. The transition from
homeomorphism rigidity to coarse rigidity, on the other hand, is mainly a matter
of elementary logic. See also \cite{VW} and \cite{KT}.
\begin{Num}\textbf{The language of Euclidean buildings\ }
\label{LanguageOfEuclideanBuildings}
Recall from \ref{LanguagesAndStructures} our language of metric space
\[
\mathcal L_\ms=\{+,\cdot,\leq,0,1,d,X,R\}.
\]
We now enlarge this
language so that we can state results about Euclidean buildings.
We fix a spherical Coxeter group $W=\bra{I\mid (ij)^{m_{ij}}=1}$
and add constants $i\in I$ for its generators and the $m_{i,j}$
to the language.
We also add a unary predicate $W$ for the Weyl group.
This allows us to include the spherical Weyl group $W$ as well as the
affine Weyl group $W\RR^n$ into our structure. 
Then we add function symbols for the coroots $\alpha_i:\RR^n\rTo\RR$,
whose kernels are the reflection hyperplanes corresponding to the
generators $i\in I$.
This allows us to describe Weyl chambers, half spaces, walls and
Weyl simplices in $\RR^n$ in our language. The last ingredient that is
missing are the coordinate charts $\phi:\RR^n\rTo X$.
The number of charts will depend very much
on the Euclidean building $X$. In order to allow some flexibility,
we view the charts as a family of maps 
\[
\{\phi_\ell:\RR^n\rTo X\}_{\ell\in L}.
\]
Thus we add a unary predicate $L$ for the indices $\ell$ that label the charts,
and one single $n+1$-ary function symbol $\phi:L\times\RR^n\rTo X$, whose first
entry is the label $\ell\in L$.
This is the language $\mathcal L_{\mathrm{eb}}$ of Euclidean buildings.

It is straightforward to see that a Euclidean building $(X,\A)$ can be viewed 
as an $\mathcal L_{\mathrm{eb}}$-structure. The axioms (EB1)--(EB5) can
be stated without difficulties in this language, 
and so can be most of the results about Euclidean buildings that we proved in 
Section~\ref{EuclideanBuildingsSection}.
Below we add a few more elements to the language $\mathcal L_{\mathrm{eb}}$,
for example the vector distance function, or a second Euclidean
building $Y$ and a coarse map $f:X\rTo Y$ between the two. We note also that
the axioms (EB5') and (EB5$^+$) can be stated if we add one function symbol 
$\rho$ for the retractions,
\[
\rho_\ell:X\rTo\phi_\ell(\RR^n)
\]
having as its first argument the index $\ell$ describing the target apartment
$A=\phi_\ell(\RR^n)$. The center of the retraction is $\phi_\ell(0)$.
\end{Num}
\begin{Num}\textbf{Ultraproducts of Euclidean buildings\ }
\label{UltraProductsOfEuclideanBuildings}
Suppose now that $K$ is a countably infinite set and that $D$ is a nonprincipal
ultrafilter on $K$. Suppose that $(X_k,\A_k)_{k\in K}$
is a family of Euclidean buildings, all of which are modeled on the same spherical
Coxeter group $W$. 
We may view them as $\mathcal L_{\mathrm{eb}}$-structures $\mathfrak X_k$
and take their ultraproduct $\mathfrak X_D$. 
We note that $W_D=W$ by \ref{UltrapowerOfFiniteStructure}.
By \L os' Theorem \ref{LosThm},
the ultraproduct satisfies the axioms (EB1)--(EB5), except that the
real numbers are now replaced by the real closed field $^*\RR$ and that
the apartments are modeled on $^*\RR^n$. The coordinate changes are now described
by the nonstandard affine Weyl group $W\ltimes {}^*\RR^n$.
One can show that $\mathfrak X_D$ is a
\emph{generalized affine building} in the sense of Bennett \cite{Bennett}
\cite{Hitz} \cite{SchwerStruyve}.
We call $\mathfrak X_D$ a 
\emph{nonstandard Euclidean building}, with \emph{nonstandard charts},
\emph{nonstandard apartments}, \emph{nonstandard Weyl simplices}, and so on.
\end{Num}
\begin{Num}\textbf{Ultralimits of Euclidean buildings\ }
\label{UltralimitsOfEuclideanBuildings}
Suppose that $\mathfrak X_D$ is an ultraproduct of of a family of
Euclidean buildings $(X_k,\A_k)$ as in \ref{LanguageOfEuclideanBuildings} and \ref{UltraProductsOfEuclideanBuildings}.
Let $p\in X_D$ and let $r>0$ be a nonstandard real. 
Then we have the ultralimit $C(X_D,p,r)$ of the underlying metric spaces $X_k$ as in
\ref{UltraLimitsAndAsymptoticCones}. 
Recall that 
\[
 X_D^{(p,r)}=\{x\in X_D\mid {\textstyle\frac1r}d(p,x)\in{}^*\RR_\fin\}
\]
and put
\[
({}^*\RR_\fin^n)^{(0,r)}=\{v\in{}^*\RR^n\mid {\textstyle\frac1r}||v||\in{}^*\RR_\fin\}=\{vr\mid v\in{}^*\RR_\fin^n\}.
\]
Suppose that $\phi_\ell$
is a nonstandard chart with $\phi_\ell(0)\in X_D^{(p,r)}$.
Then we have a commuting diagram
\begin{diagram}[height=2em]
 ({}^*\RR_\fin^n)^{(0,r)}&&\rTo^{\phi_\ell}&&X_D^{(p,r)}\\
 \dTo&&&&\dTo\\
 \RR^n&&\rDotsto^{\tilde\phi_\ell}&&C(X_D,p,r)
\end{diagram}
where the vertical arrows identify points at infinitesimal $\frac1rd$-distance.
The map $\tilde\phi_\ell$ is an isometric injection, and we let
$\tilde\A$ denote the set of all maps $\tilde\phi_\ell$ that arise
in this way. The corresponding affine apartments (the images of the $\tilde\phi_\ell$)
will be denoted by
$\tilde A\subseteq C(X_D,p,r)$.
\end{Num}
The following result is proved in \cite[5.1.1]{KL}. A more general result
is proved by Schwer and Struyve in \cite[6.1]{SchwerStruyve}.
We remark that $X_D^{(p,r)}$ is also a generalized affine building in the sense of Bennett \cite{Bennett}.
\begin{Thm}
\label{UltralimitsOfEuclideanBuildingsThm}
Let $(X_k,\A_k)_{k\in K}$ be a countably infinite family of Euclidean buildings,
with a fixed spherical Coxeter group $W$.
Let $C(X_D,p,r)$ and $\tilde\A$ be as in \ref{UltralimitsOfEuclideanBuildings}.
Then the ultralimit $C(X_D,p,r)$ is a metrically complete Euclidean building,
with $\tilde\A$ as its complete set of apartments.

\proof
From the construction of $\tilde\A$ and \L os' Theorem it is clear that
$\tilde\A$ satisfies the axioms (EB1) and (EB2). The retractions $X_D\rTo A$
onto the nonstandard affine apartments are $1$-Lipschitz and descend therefore 
to retractions $C(X_D,p,r)\rTo\tilde A$, hence axiom (EB5$^+$) holds.
The remaining two axioms require more work.

\medskip\noindent
\emph{Axiom (EB3) holds.}
We note that  $C(X_D,p,r)$ is CAT$(0)$ by
\ref{UltralimitsPreserveCurvature}. Since the affine apartments are convex,
their intersections are also convex. Now we add the vector distance $\Theta$ to
the structure. It is clear that this gives us a nonstandard vector distance 
on $X_D^{(p,r)}\subseteq X_D$. Suppose that $x,x',y\in X_D^{(p,r)}$
and that $x$ and $x'$ have infinitesimal $\frac1rd$-distance. Then
$\Theta(y,x)$ and $\Theta(y,x')$ also have infinitesimal $\frac1rd$-distance. 
From the symmetry of the vector distance \ref{VectorDistance},
we conclude that $\Theta(x,y)$ and $\Theta(x',y)$ have infinitesimal distance.
This implies that the nonstandard vector distance descends to an ordinary vector
distance $\tilde\Theta:C(X_D,p,r)\times C(X_D,p,r)\rTo a_0$. It now follows from
\ref{LocalCompatibilityCriterion} that coordinate changes between charts in
$C(X_D,p,r)$ are given by elements of $W\RR^n$, hence (EB3) holds.

\medskip\noindent
\emph{Axiom (EB4) holds.}
We first show two auxiliary results.

\smallskip\noindent
\emph{Claim 1. Let $c,c'\subseteq X_D$ be nonstandard Weyl chambers with
$\partial c=\partial c'$. If the tips $z$ and $z'$ of 
$c$ and $c'$ are in $X_D^{(p,r)}$,
then there exists a nonstandard Weyl chamber $a\subseteq c\cap c'$
whose tip $w$ is also in $X_D^{(p,r)}$.}

\smallskip\noindent
Let $[z,\zeta_c)$ be the nonstandard central ray of $c$ (we may add these
rays to our structure and language). This ray intersects the Weyl chamber
$c'$ and there
is a unique point\footnote{The point exists by \L os' Theorem because it is 
\emph{definable} in the
structure, hence it does not matter that the ordered field
$^*\RR$ is non-archimedean.} 
$w$ with $[z,\zeta_c)\cap c'=[w,\zeta_c)$. Let
\[
s={\textstyle\frac1r}d(w,z)\in{}^*\RR.
\]
We claim that $s$ is a finite nonstandard real.
Suppose to the contrary that $s$ is infinite. Let $[x,\zeta_c)\subseteq c'$
be the maximal ray in $c'$ extending $[w,\zeta_c)$. Then the three points $z,z',x$
have mutually infinitesimal $\frac1{rs}d$-distance, while
$w$ and $z$ have $\frac1{rs}d$-distance $1 \pmod{{}^*\RR_\infi}$.
It follows that the three 
nonstandard Alexandrov angles at $w$ between the segments $[w,z],[w,z'],[w,x]$
are infinitesimally small (we may also add the Alexandrov angles to our structure
and language, so they are defined in $X_D$). 
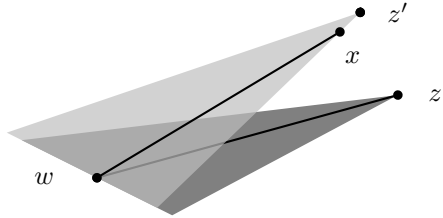
\begin{center}
\begin{pspicture}(4,0)(4,3)
\pspolygon[linestyle=none,fillstyle=solid,fillcolor=gray](3,0)(1,1)(6,1.6)
\psline[linestyle=solid]{*-*}(2,0.5)(6,1.6)
\pspolygon[linestyle=none,fillstyle=solid,fillcolor=lightgray,opacity=0.7](2.8,0.1)(0.8,1.1)(5.5,2.7)
\psline[linestyle=solid]{*-*}(5.5,2.7)(5.5,2.7)
\psline[linestyle=solid]{*-*}(2,0.5)(5.23,2.44)
\rput(1.3,0.5){\small$w$}
\rput(6.5,1.6){\small$z$}
\rput(6,2.7){\small$z'$}
\rput(5.4,2.1){\small$x$}
\end{pspicture}
\end{center}
This contradicts \ref{AngelRigidity}.
Thus $s$ is finite and hence $w\in X_D^{(p,r)}$. Now let $a$ be the unique
nonstandard $w$-based Weyl chamber with $\partial a=\partial c$.
Thus Claim 1 is proved.

\smallskip
We note that every Weyl chamber $\tilde c\subseteq C(X_D,p,r)$ arises from
some nonstandard Weyl chamber $c$ whose tip is in $X_D^{(p,r)}$.
Suppose that $\tilde c,\tilde c'\subseteq C(X_D,p,r)$ are Weyl chambers.
We choose corresponding nonstandard Weyl chambers $c,c'\subseteq X_D$ in such
a way that $\partial c,\partial c'$ have minimal gallery distance $m\in\mathbb N$
in the spherical building at infinity of $X_D$. 

\smallskip\noindent
\emph{Claim 2. In this situation, there is an nonstandard affine apartment
$A$ containing a point $w\in X_D^{(p,r)}$ and a gallery of $w$-based
nonstandard Weyl chambers $c_0,\ldots,c_m$, with $\partial c=\partial c_0$ and
$\partial c'=\partial c_m$.}

\smallskip\noindent
We proceed by induction on $m$, the case $m=0$ being trivial.
Assume now that $c_0,\ldots,c_{m-1}$ are $w$-based Weyl chambers 
contained in a nonstandard affine apartment $A$, with $w\in X_D^{(p,r)}$.
Let $a$ be the unique $w$-based Weyl chamber with $\partial a=\partial c'$.
If $a$ is contained in $A$ we are done, with $c_m=a$.
If $a$ is not contained in $A$, then we let $M\subseteq A$ denote the nonstandard
wall that separates $a$ and $c_{m-1}$. Let $H\subseteq A$
denote the closed half space bounded by $M$ whose boundary at infinity contains
$\partial c_{m-1}$. If $H$ and $X_D^{(p,r)}$ have a
point $w'$ in common, then we replace $c_0,\ldots,c_{m-1}$ by their 
$w'$-based translates in $A$, and we let $c_m$ denote the unique 
$w'$-based Weyl chamber with $\partial c_m=\partial c'$.
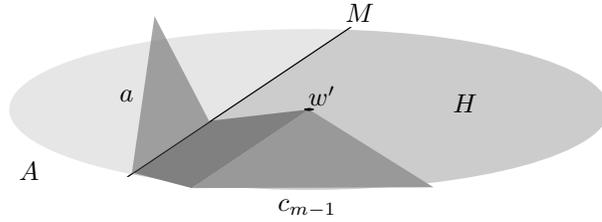
\begin{center}
\begin{pspicture}(-2,-1.3)(4,2)
\psset{viewpoint=-5 -2 1.5}
\ThreeDput[normal=0 0 1]{
\psclip{\psframe[linestyle=none](-5,-5)(5,1)}
\pscircle[linestyle=none,fillstyle=solid,fillcolor=gray8](0,0){4}
\endpsclip
\psclip{\psframe[linestyle=none](-5,1)(5,5)}
\pscircle[linestyle=none,fillstyle=solid,fillcolor=gray9](0,0){4}
\endpsclip
\pspolygon[linestyle=none,fillstyle=solid,fillcolor=gray6](0,0)(-4.2,0)(-2.969,-2.969)
\pspolygon[linestyle=none,fillstyle=solid,fillcolor=gray](0,0)(-4.2,0)(-3.785,1)(-1,1)
\psline[linestyle=solid](-4,1)(4,1)
\psline[linestyle=solid]{*-*}(0,0)(0,0)
}
\ThreeDput[normal=0 -1 0](0,1,0){
\pspolygon[linestyle=none,fillstyle=solid,fillcolor=gray,opacity=0.7](-3.785,0)(-1,0)(-2.969,1.969)
}
\rput(.2,.2){\small$w'$}
\rput(0.7,1.3){\small$M$}
\rput(2.1,0.1){\small$H$}
\rput(-3.7,-.8){\small$A$}
\rput(0,-1.3){\small$c_{m-1}$}
\rput(-2.4,0.2){\small$a$}
\end{pspicture}
\end{center}
Let $\bar a$ be the unique $w'$-based Weyl chamber in $H$ whose germ near $w'$ is opposite 
to the germ of $a$. Then $a,\bar a$ are contained in a unique nonstandard affine
apartment $A'$ containing $c_0\cup\cdots\cup c_m$, hence we are done.
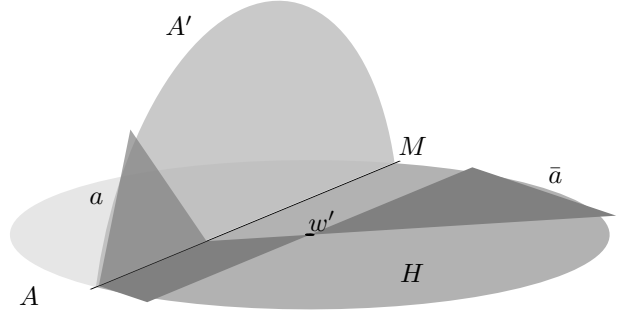
\begin{center}
\begin{pspicture}(-2,-1)(4,3)
\psset{viewpoint=-5 -3 1.5}
\ThreeDput[normal=0 0 1]{
\psclip{\psframe[linestyle=none](-5,-5)(5,1)}
\pscircle[linestyle=none,fillstyle=solid,fillcolor=gray7](0,0){4}
\endpsclip
\psclip{\psframe[linestyle=none](-5,1)(5,5)}
\pscircle[linestyle=none,fillstyle=solid,fillcolor=gray9](0,0){4}
\endpsclip
\pspolygon[linestyle=none,fillstyle=solid,fillcolor=gray](0,0)(4.2,0)(2.969,-2.969)
\pspolygon[linestyle=none,fillstyle=solid,fillcolor=gray](0,0)(-4.2,0)(-3.785,1)(-1,1)
\psline[linestyle=solid](-4,1)(4,1)
\psline[linestyle=solid]{*-*}(0,0)(0,0)
}
\ThreeDput[normal=0 -1 0](0,1,0){
\psclip{\psframe[linestyle=none](-5,0)(5,5)}
\pscircle[linestyle=none,fillstyle=solid,fillcolor=gray7,opacity=.7](0,-1){4}
\endpsclip
\pspolygon[linestyle=none,fillstyle=solid,fillcolor=gray,opacity=0.7](-3.785,0)(-1,0)(-2.969,1.969)
}
\rput(.2,.2){\small$w'$}
\rput(1.4,1.2){\small$M$}
\rput(1.4,-0.5){\small$H$}
\rput(-3.7,-.8){\small$A$}
\rput(3.3,.8){\small$\bar a$}
\rput(-2.8,.5){\small$a$}
\rput(-1.7,2.8){\small$A'$}
\end{pspicture}
\end{center}
Suppose that $H\cap X_D^{(p,r)}=\emptyset$. Then $M$ separates
$w$ from $\partial c_{m-1}$, hence $a$ has the same germ near $w$ as
$c_{m-1}$. 
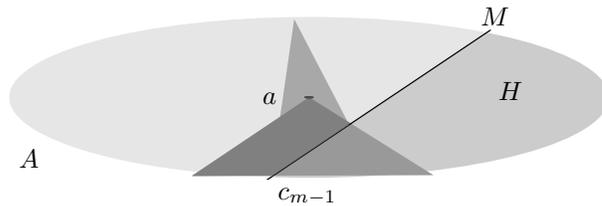
\begin{center}
\begin{pspicture}(-2,-1.5)(4,1)
\psset{viewpoint=-5 -2 1.5}
\ThreeDput[normal=0 0 1]{
\psclip{\psframe[linestyle=none](-5,-5)(5,-1)}
\pscircle[linestyle=none,fillstyle=solid,fillcolor=gray8](0,0){4}
\endpsclip
\psclip{\psframe[linestyle=none](-5,-1)(5,5)}
\pscircle[linestyle=none,fillstyle=solid,fillcolor=gray9](0,0){4}
\endpsclip
\pspolygon[linestyle=none,fillstyle=solid,fillcolor=gray6](0,0)(-4.2,0)(-2.969,-2.969)
\pspolygon[linestyle=none,fillstyle=solid,fillcolor=gray](0,0)(-4.2,0)(-3.785,-1)(-1,-1)
\psline[linestyle=solid](-4,-1)(4,-1)
\psline[linestyle=solid]{*-*}(0,0)(0,0)
}
\ThreeDput[normal=0 -1 0](0,-1,0){
\pspolygon[linestyle=none,fillstyle=solid,fillcolor=gray,opacity=0.6](-3.785,0)(-1,0)(-2.969,1.969)
}
\rput(-.5,0){\small$a$}
\rput(2.5,1.1){\small$M$}
\rput(2.7,0.1){\small$H$}
\rput(-3.7,-.8){\small$A$}
\rput(0,-1.3){\small$c_{m-1}$}
\end{pspicture}
\end{center}
This situation is excluded by our assumptions.
Thus we have proved Claim 2.

\smallskip
In the setting of Claim 2, it follows from Claim 1 that $\tilde c\cap \tilde c_0$
contains a Weyl chamber, and that $\tilde c'\cap \tilde c_m$ contains a
Weyl chamber. Therefore (EB4) holds.

\medskip\noindent
\emph{The atlas $\tilde\A$ is maximal.}
We have to show that every geodesic line and ray is contained
in some affine apartment \cite[2.18]{Par}.
Let $\gamma:\RR\rTo C(X_D,p,r)$ be a geodesic line. We choose points
$x_k\in X_D^{(p,r)}$ corresponding to $\gamma(k)\in C(X_D,p,r)$, for $k=0,\pm1,\pm2,\ldots$.
For $k\geq 0$, let $A_k$ be a nonstandard affine apartment containing $\{x_{-k},x_k\}$.
Then $x_\ell$ is $\frac1rd$-infinitesimally close to $A_k$ for all $\ell=-k,\ldots,k$.
Hence we have the following property of the set $T=\{x_k\mid k\in\ZZ\}$.
For every finite subset $S\subseteq T$ there exists a nonstandard affine apartment 
$A_S$ such that all elements of $S$ are $\frac1rd$-infinitesimally close
to $A_S$. By \ref{overspill} there exists a nonstandard affine apartment
$A_T$ such that all members of $T$ are $\frac1rd$-infinitesimally close to $A_T$.
Thus $\gamma(k)\in\tilde A_T$ holds for all $k\in\ZZ$. Since $\tilde A_T$ is convex,
we have $\gamma(\RR)\subseteq\tilde A_T$. The reasoning for rays is completely
analogous.
\qed
\end{Thm}
The following consequence is often useful. Note that the metric completion of
a non-complete Euclidean building need not be a Euclidean building
\cite[6.9]{KramerLocal},
see also \cite{MSSS}.
\begin{Cor}
Let $(X,\A)$ be a Euclidean building. Then there exists a complete Euclidean
building $(\bar X,\bar\A)$ of the same type, with $X\subseteq\bar X$ and
$\A\subseteq \bar\A$. The action of the automorphism group of $(X,\A)$ 
extends to an action on $(\bar X,\bar\A)$.

\proof
Let $X_D$ be the ultrapower of $X$ and put $\bar X=C(X_D,p,r)$. The diagonal
embedding $X\rTo X_D$ extends to an embedding $X\rTo\bar X$. 
The automorphism group of $(X,\A)$ acts diagonally on $X_D$, on $X_D^{(p,r)}$,
and hence on $\bar X$ in a natural way.
\qed
\end{Cor}

\begin{Num}\textbf{The structure \boldmath$\mathfrak C$\ }
\label{TheStructureC}
Suppose that $f:X\rTo Y$ is a coarse equivalence of
Euclidean buildings. We consider the structure $\mathfrak C$ 
consisting of the two Euclidean buildings and the map $f$,
and we take the ultrapower $\mathfrak C_D$ of this structure.
If $p\in X_D$ is any basepoint and if $r\in{}^\ast\RR$ is 
infinitely large, then $f$ induces a bi-Lipschitz homeomorphisms
\[
C(f):C(X_D,p,r)\rTo C(Y_D,f(p),r)
\] 
by \ref{CoarseToLipschitz}. By \ref{UltralimitsOfEuclideanBuildingsThm},
the asymptotic cones  $C(X_D,p,r)$ and $C(Y_D,f(p),r)$ are metrically
complete Euclidean buildings with complete apartment systems, and 
$C(f)$ maps affine apartments to affine apartments by
\ref{TopologicalRigidityOfApartments}.
We also note the following.
The ordered field $^\ast\RR$ is non-archimedean, so a bounded
set of nonstandard reals will in general not have a supremum.
Nevertheless, the quantity 
\[
Hd_{B_r(y)}(f(A),A')=Hd(f(A)\cap B_r(y),A'\cap B_r(y))
\]
is defined for all 
nonstandard affine apartments $A\subseteq X_D$,
$A'\subseteq Y_D$, points $y\in Y_D$ and positive nonstandard reals $r>0$.
The reason for this is that we can add $Hd_{B_r(y)}(f(A),A')$
as an $\RR\cup\{\infty\}$-valued function (depending on arguments
$\ell,\ell',y,r$, with $A=\phi_\ell(\RR^n)$, $A'=\phi'_{\ell'}(\RR^n)$) 
to the structure $\mathfrak C$ and then take its 
ultrapower. \L os' Theorem \ref{LosThm} guarantees that this function
has exactly the intended meaning in $\mathfrak C_D$, namely that
of an $^\ast\RR$-valued Hausdorff distance%
\footnote{The point is that
we measure only the Hausdorff distance between so-called 
\emph{internal} or \emph{definable} sets.}.
These observations are the main steps in our proof of the following theorem
\cite[7.1.5]{KL}.
\end{Num}
\begin{Thm}[Higher dimensional Morse Lemma]
\label{HigherDimensionalMorseLemma}
Let 
\[
f:X\rTo Y
\]
be a coarse equivalence between Euclidean buildings. There
exists a real constant $r_f>0$ such that the following holds. For every affine apartment
$A\subseteq X$ there is a unique affine apartment $A'$ in the complete 
apartment system of $Y$ with 
\[
Hd(f(A),A')\leq r_f.
\]
\end{Thm}
We first prove a weaker statement in the ultraproduct.
\begin{Lem}
\label{Lemma1Morse}
Let $\mathfrak C$ be as in \ref{TheStructureC}. Let
$A\subseteq X_D$ be a nonstandard affine apartment and let $y \in f(A)$.
Let $r\geq 1$ be a nonstandard real. 
There exist a nonstandard affine apartment $A'\subseteq Y_D$
and a finite nonstandard real number $s\geq 0$ such that 
\[
Hd_{B_{\eps r}(y)}(f_D(A),A')\leq s
\]
holds for all nonstandard $\eps$ with $\frac12\leq\eps\leq1$.

\proof
This is clear if $r$ is finite: we may choose any nonstandard
affine apartment $A'$ containing $y$ and choose any real number $s\geq 2r$.  So suppose that
$r$ is infinite. We choose a preimage $x\in A$ of $y$. 
Since $r$ is infinite, we have by Lemma
\ref{CoarseToLipschitz} a bi-Lipschitz map
\[
C(f):C(X_D,x,r)\rTo C(Y_D,y,r).
\]
Let $\tilde A\subseteq C(X_D,x,r)$ denote the affine corresponding to $A$. 
By Theorem \ref{TopologicalRigidityOfApartments}
$C(f)$ maps $\tilde A$ onto an affine apartment $\tilde A'$
in the complete apartment system of the Euclidean
building $C(Y_D,y,r)$. By \ref{UltralimitsOfEuclideanBuildingsThm} 
there is a nonstandard affine apartment $A'\subseteq Y_D$ corresponding
to $\tilde A'$. Suppose now that $\eps$ is a nonstandard real with $\frac12\leq\eps\leq 1$.
Let 
\[
s_\eps=Hd_{B_{\eps r}(y)}(f(A),A').
\]
Since we have in $C(Y_D,y,r)$ that
\[
B_\eps(\tilde y)\cap C(f)(\tilde A)=B_\eps(\tilde y)\cap \tilde A',
\]
the quotient $s_\eps/r$ is infinitesimally small.
We claim that $s_\eps$ is finite. Suppose to the contrary that $s_\eps>0$ is infinite. 
There is either a point $z\in A'\cap B_{\eps r}(y)$ such that 
$B_{s_\eps/2}(z)\cap f(A)=\emptyset$ or a point $z\in f(A)\cap B_{\eps r}(y)$ such that
$B_{s_\eps/2}(z)\cap A'=\emptyset$. Because $s_\eps$ is infinite, we have a bi-Lipschitz map
\[
C(f):C(X_D,w,s_\eps)\rTo C(Y_D,z,s_\eps)
\]
(for a suitably chosen point $w\in X_D$).  In $C(Y_D,z,s_\eps)$, the sets
$C(f)(\tilde A)$ and $\tilde A'$ are affine apartments, and $C(f)(\tilde A)\neq\tilde A'$
by the choice of $z$. It is also clear from the construction that 
$C(f)(\tilde A)$ and $\tilde A'$ have Hausdorff distance at most $1$ in $C(Y_D,z,s_\eps)$.
This is impossible by \cite[4.6.4]{KL} \cite[p.~10]{Par}, or
by \ref{CloseApartmentsIntersect}. Thus $s_\eps$ has to be finite.

Now we claim that the set $\{s_\eps\mid\frac12\leq\eps\leq 1\}\subseteq{}^*\RR$ 
has a finite upper bound.
Suppose that this is false. Then we find for every $k\in\NN$ an $\eps_k$ 
with $\frac12\leq\eps_k\leq 1$ such that $s_{\eps_k}\geq k$. 
By \ref{overspill} we find an $\eps$ with $\frac12\leq\eps\leq1$ such that $s_\eps$
is infinite, a contradiction. Hence $\{s_\eps\mid\frac12\leq\eps\leq 1\}\subseteq{}^*\RR$ 
has a finite upper bound $s$.
\qed
\end{Lem}
The next lemma says that the numbers $s$ occurring in the previous lemma
can be bounded uniformly from above. This is another application of \ref{overspill}.
\begin{Lem}
\label{Lemma2Morse}
Let $\mathfrak C$ be as in \ref{TheStructureC}. 
There exists a finite real constant $s\geq 0$ such that the following holds. 
For every $r\geq1$,
every nonstandard affine apartment $A\subseteq X_D$ and every $x\in A$ there exists
a nonstandard affine apartment $A'\subseteq Y_D$ such that 
$Hd_{B_{\eps r}(f(x))}(f(A),A')\leq s$ holds for all $\frac12\leq\eps\leq 1$.

\proof
Assume that this is false. Then we can find for every $k\in\NN$ a positive nonstandard real
$r_k\geq 1$, an $\eps_k$ with $\frac12\leq\eps_k\leq 1$, a nonstandard affine apartment
$A_k\subseteq X_D$, a point $x_k\in A_k$ such that for every nonstandard affine apartment 
$A'\subseteq Y_D$ we have $Hd_{B_{\eps_k r_k}(f(x_k))}(f(A_k),A')\geq k$.
By \ref{overspill} there exist a nonstandard affine apartment $A\subseteq X_D$,
a point $x\in A$,
a nonstandard real $r\geq 1$ and an $\eps$ with $\frac12\leq\eps\leq 1$ such that 
$Hd_{B_{\eps r}(f(x))}(f(A),A')$ is infinite for all nonstandard affine apartments $A'\subseteq Y_D$.
This contradicts Lemma \ref{Lemma1Morse}.
\qed
\end{Lem}
By \L os' Theorem \ref{LosThm} we have the following immediate consequence for the
coarse equivalence that we started with.
\begin{Cor}
\label{CorMorse3}
Let 
\[
f:X\rTo Y
\]
be a a coarse equivalence between Euclidean buildings.
There exists a real constant $r_f\geq 0$ such that the following holds. 
For every $r\geq 1$,
every affine apartment $A\subseteq X$, every $x\in A$ there exists
an affine apartment $A'\subseteq Y$ such that 
\[
Hd_{B_{\eps r}(f(x))}(f(A),A')\leq r_f
\]
holds for all $\frac12\leq\eps\leq 1$.
\qed
\end{Cor}
\emph{Proof of the higher dimensional Morse Lemma \ref{HigherDimensionalMorseLemma}.}
Let $A\subseteq X$ be an affine apartment, let $x\in A$ and put 
$Z=f(A)$ and $z=f(x)$. Let $r_f$ be as in \ref{CorMorse3}.
For every $s\geq1$ we may by \ref{CorMorse3} choose an
affine apartment $A_s\subseteq Y$ such that 
\[
Hd_{B_{\eps s}(z)}(Z,A_s)\leq r_f
\]
holds for all $s\in[\frac12,1]$.
Let $c=c_{2r_f}$ be the constant from \ref{CloseApartmentsIntersect},
corresponding to Hausdorff distance $2r_f$.
For $s\geq 1$ we have
$Hd_{B_{s+c}(z)}(Z,A_{s+c})\leq r_f$ and 
$Hd_{B_{s+c}(z)}(Z,A_{2(s+c)})\leq r_f$, whence
\[
Hd_{B_{s+c}(z)}(A_{s+c},A_{2(s+c)})\leq 2r_f
\]
and, by \ref{CloseApartmentsIntersect},
\[
 B_s(z)\cap A_{s+c}=B_s(z)\cap A_{2(s+c)}.
\]
Let $w=\pi_{A_{s+c}}(z)=\pi_{A_{2(s+c)}}(z)$. We have then, for $s>r_f$,
\[
B_{s-r_f}(w)\cap A_{s+c}=B_{s-r_f}(w)\cap A_{2(s+c)}.
\]
Then $A_{r_f+c},A_{2(r_f+c)},A_{4(r_f+c)},\ldots,A_{2^k(r_f+c)},\ldots$
gives us a nested sequence of metric open $n$-balls (where $n=dim(X)$)
with centers $w$ and of radii $2^k(r_f+c)-c$. Their union 
\[
 A'=\bigcup_{k=0}^\infty \left(B_{2^k(r_f+c)-c}(w)\cap A_{2^k(r_f+c)}\right)
\]
is isometric to $\RR^n$ and hence by \ref{MaximalFlatsProp} an affine apartment in the complete
apartment system of $Y$. From the construction of $A'$ it is clear that
\[
Hd(Z,A')\leq r_f.
\]
\qed

\medskip\noindent
We also have the following.
\begin{Prop}
\label{CoarseEqivalencesPreserveEuclideanFactors}
Let $X,Y$ be Euclidean buildings whose
spherical buildings at infinity are thick. Let
$f:X\times\RR^{m_1}\rTo Y\times\RR^{m_2}$ be a
coarse equivalence. Then $m_1=m_2$.

\proof
We choose an infinite positive nonstandard real $r$ and a thick point $p\in X$, see \ref{ThickPointsExist}.
We then consider the asymptotic cone 
\[
C((X\times\RR^{m_1})_D,r,p\times 0)\cong C(X_D,r,p)\times\RR^{m_1}
\]
with respect to the constant 
families $X_k=X$ and $p_k=p$. It is clear that the constant sequence $(p_k)_{k\in K}$ 
represents a thick point in the ultraproduct
$X_D$, and hence also in the asymptotic cone $C(X_D,r,p)$. In particular, the spherical building at infinity
of the Euclidean building $C(X_D,r,p)$ is thick. Now we may apply
Corollary \ref{HomeosPreserveEuclideanFactors} to the continuous map $C(f)$.
\qed
\end{Prop}
Now we prove that a coarse equivalence
of Euclidean buildings induces an isomorphism between the
spherical buildings at infinity. In order to show this,
we have to describe the spherical building at infinity
by coarse data. 
\begin{Lem}
\label{NbhdInts}
Let $A_1,\ldots, A_k$ be affine apartments in a Euclidean
building $X$. The following are equivalent.

(i) $\partial A_1\cap\cdots\cap\partial A_k\neq\emptyset$.

(ii) There is an unbounded set which is dominated by each of the 
affine apartments $A_1,\ldots,A_k$.

\proof
To see that (i) implies (ii), 
let $a\subseteq X$ be a Weyl simplex representing an element
$\partial a\in\partial A_1\cap\cdots\cap\partial A_k$.
For each $j$ there is a Weyl simplex $a_j\subseteq A_j$ 
representing $\partial a$. Since $a$ and $a_j$ have finite
Hausdorff distance, each $A_j$ dominates $a$.

Before we prove the converse implication, we note the following.
If an affine apartment $A$ dominates a set $Z\subseteq X$, then
$\pi_AZ$ has finite Hausdorff distance from $Z$. If $Z$ is unbounded,
then there exists a Weyl simplex $a\subseteq A$ that contains
an unbounded subset of $\pi_AZ$. In particular, there exists then a
Weyl simplex $a\subseteq A$ that dominates an unbounded subset of 
$Z$ (or $\pi_AZ$).

Now we assume that the unbounded set
$Z\subseteq X$ is dominated by the affine apartments $A_1,\ldots,A_k$. 
Consider the unbounded set $Y=\pi_{A_1}Z$.
There exists a Weyl simplex
$a\subseteq A_1$ of minimal dimension which
dominates an unbounded subset $Y_1\subseteq Y$.
We claim that 
\[
\partial a\in \partial A_1\cap\cdots\cap\partial A_k.
\]
Let $j>1$. In $A_j$ we find a Weyl chamber $c_j$
which dominates an unbounded subset
$Y_j$ of $\pi_{A_j}(Y_1)$. Let $A_j'$ be an affine
apartment containing representatives $a'$ and $c_j'$ of 
$\partial a$ and $\partial c_j$. Then
$Y_j'=\pi_{A_j'}(Y_j)$ has finite Hausdorff distance from
$Y_j$.
Since both $a'$ and $c_j'$ dominate the unbounded set
$Y_j'$, and since $a'$ is also a Weyl simplex of minimal dimension
which dominates an unbounded subset of $Y$,
we have $\partial a\subseteq\partial c_j\in\partial A_j$.
\qed
\end{Lem}
Suppose that $X$ and $Y$ are Euclidean buildings and
that $f:X\times\RR^{m_1}\rTo Y\times\RR^{m_2}$ is
a coarse equivalence. We may view $\RR^{m_1}$ as a Euclidean
building. Then every affine apartment in $X\times\RR^{m_1}$
is of the form $A\times\RR^{m_1}$, where
$A\subseteq X$ is an affine apartment.
By \ref{HigherDimensionalMorseLemma} there is
a map $f_*$ from the set of all affine apartments of $X$ to the 
set of all affine
apartments of $Y$, such that $f(A\times\RR^{m_1})$ has finite
Hausdorff distance from $f_*A\times \RR^{m_2}$.
If $g$ is a coarse inverse of $f$, then $g_*$ is an inverse
of $f_*$.
\begin{Lem}
\label{AptCplxHom}
Suppose that $X$ and $Y$ are Euclidean buildings. Let 
\[
f:X\times\RR^m\rTo\allowbreak Y\times\RR^m
\]
be a
coarse equivalence. Then $f_*$ induces an isomorphism
between the apartment complexes
$AC(\partialmax X)$ and $AC(\partialmax Y)$.

\proof
Let $A_1,\ldots,A_k\subseteq X$ be affine apartments with
$\partial A_1\cap\cdots\cap \partial A_k\neq\emptyset$.
Let $a$ be a $1$-dimensional Weyl simplex representing 
a vertex $\partial a\in \partial A_1\cap\cdots\cap \partial A_k$.
Choose $s>0$ such that $a\subseteq B_s(A_1)\cap \cdots\cap B_s(A_k)$.
Let $r_f>0$ be as in \ref{HigherDimensionalMorseLemma} and let
\[
Z=B_{r_f+\rho(s)}(f_*A_1)\cap\cdots\cap B_{r_f+\rho(s)}(f_*A_k)\subseteq Y,
\]
where $\rho$ is the control function for $f$.
Then $f$ restricts to a quasi-isometric embedding 
of $a\times\RR^{m}$ into $Z\times\RR^{m}$.
If $Z$ is bounded, we get a quasi-isometric
embedding of $a\times\RR^{m}\cong[0,\infty)\times\RR^m$ into
$\RR^{m}$, which is impossible
(by topological dimension invariance, applied to the
asymptotic cones, see \ref{DimensionInvariance} and the proof of 
\ref{BoundaryLemma}). Hence $Z$ is unbounded.
Therefore $\{f_*A_1,\ldots,f_*A_k\}$ is a simplex in
$AC(\partialmax X_2)$ by \ref{NbhdInts}. 
It follows that $f_*$ is a simplicial
map $f_*:AC(\partialmax X_1)\rTo AC(\partialmax X_2)$.
If $g$ is a coarse inverse of $f$, then $g_*$ is a simplicial
inverse of $f_*$.
\qed
\end{Lem}
Combining these results, we have the following first main result
about coarse equivalences between Euclidean buildings. This is
Theorem~II in the introduction.
\begin{Thm}
\label{MainThmII}
Let $X_1$ and $X_2$ be Euclidean buildings whose spherical
buildings at infinity $\partialmax X_1$ and
$\partialmax X_2$ are thick. Let
\[
f:X_1\times\RR^{m_1}\rTo X_2\times\RR^{m_2}
\]
be a coarse
equivalence. Then $m_1=m_2$ and the map $f_*$ on the
affine apartments extends uniquely to a simplicial isomorphism
$f_*:\partialmax X_1\rTo\partialmax X_2$.

\proof
By \ref{CoarseEqivalencesPreserveEuclideanFactors} we have $m_1=m_2$.
By \ref{AptCplxHom} the map $f_*$ induces a simplicial isomorphism
between $AC(\partialmax X_1)$ and $AC(\partialmax X_2)$.
By \ref{AptCplxProp}
$f_*$ induces a simplicial isomorphism
$f_*:\partialmax X_1\rTo^\cong \partialmax X_2$.
\qed
\end{Thm}
The thickness of the spherical buildings is essential for the
argument. However, the Euclidean factors which are allowed
in the theorem lead also to a result for the case that the
spherical buildings at infinity are weak buildings.
\begin{Cor}
Let $f:X_1\rTo X_2$ be a coarse equivalence between Euclidean buildings.
Then the induced map on the affine apartments induces an isomorphism
between the thick building factors in the
reductions of $\partialmax X_1$ and $\partialmax X_2$.

\proof
This follows from \ref{ThickReductionEuclideanBuilding},
and \ref{MainThmII}.
\qed
\end{Cor}
Note that we do not claim that $f$ induces directly a map 
$\partial f:\partial X_1\rTo\partial X_2$
between the Tits boundaries  in the sense of CAT$(0)$ geometry.
This will in general not be the case; for example, $f$ could be a bi-Lipschitz
homeomorphism of the Euclidean cone $\EB (B)$ over a thick spherical building $B$.
Such a self homeomorphism can be rather wild at infinity. Our combinatorial
construction of $f_*$ applies nevertheless.
\begin{Num}
\label{SetupForThmIII}
It remains to prove Theorem III. We fix some notation.
We assume that $X_1$ and $X_2$ are metrically complete
Euclidean buildings whose spherical buildings at infinity 
$\partialmax X_1$ and $\partialmax X_2$ are thick.
Furthermore, we assume that 
\[
f:X_1\times\RR^{m_1}\rTo X_2\times\RR^{m_2}
\]
is a coarse equivalence. By \ref{MainThmII}, $f$ 
induces an isomorphism 
\[
f_*:\partialmax X_1\rTo_\cong\partialmax X_2
\]
which is characterized by the fact that for every affine
apartment $A\subseteq X_1$, the image $f(A\times \RR^{m_1})$
has finite Hausdorff distance from $f_*A\times\RR^{m_2}$.
We put 
\[
n=\dim(X_1)=\dim(X_2).
\]
The Euclidean factors have by \ref{CoarseEqivalencesPreserveEuclideanFactors}
the same dimension, which we denote by 
\[
m=m_1=m_2. 
\]
Finally, we put
\[
f(x\times v)=f_1(x\times v)\times f_2(x\times v).
\]
Let $a,b\in\partialmax X_1$ be opposite panels. As in
\ref{WallTrees}, we denote by
$X_1(a,b)$ the union of all affine apartments in
$X_1$ whose
boundary contains $a$ and $b$. Consider the closed convex subsets
\[
Y_1=X_1(a,b)\times\RR^m\subseteq X_1\times\RR^m\quad\text{ and }\quad
Y_2=X_2(f_*a,f_*b)\times\RR^m\subseteq X_2\times\RR^m.
\]
The composite $\pi_{Y_2}\circ f:Y_1\rTo Y_2$ is a controlled map.
If $A\subseteq X_1(a,b)$ is an affine apartment, then
$f_*A$ is an affine apartment in $X_2(f_*a,f_*b)$;
see \ref{AptCplxHom}. By
\ref{HigherDimensionalMorseLemma}, there is a uniform constant $r_f>0$
such $d(f(y),\pi_{Y_2}(f(y))\leq r_f$ for all $y\in Y_1$,
so $f|Y_1$ and $\pi_{Y_2}\circ f|Y_1$ have finite distance.
If $g$ is a coarse inverse of $f$, then
$\pi_{Y_1}\circ g|_{Y_1}$ is therefore a coarse inverse of 
$\pi_{Y_2}\circ f|_{Y_2}$,
and we obtain a coarse equivalence 
\[
\pi_{Y_2}\circ f:Y_1\rTo Y_2.
\]
By \ref{WallTrees}, $Y_i$ factors for $i=1,2$
as a metric product $T_i\times\RR^{n-1+m}$
of a metrically complete leafless tree $T_i$, the wall tree, and a Euclidean space.

The ends of the wall tree $T_1$ correspond bijectively
to the chambers of $\partialmax X_1$ containing the panel $a$.
If the type of the panel $a$ is not isolated in the Coxeter
diagram of $\partialmax X_1$, then $\Pi^+_{\partialmax X_1}(a)$
acts $2$-transitively on the ends of $T_1$.
The boundary map
$f_*:\partialmax X_1\rTo\partialmax X_2$ is clearly
equivariant with respect to the isomorphism
$\Pi^+_{\partialmax X_1}(a)\rTo\allowbreak\Pi^+_{\partialmax X_2}(f_*a)$.
In this situation we may apply \ref{MainTreeProp} and
conclude that $f_*$ extends to  an isometry from $T_1$ to
$T_2$, possibly after rescaling.
\end{Num}
\begin{Lem}
\label{IrredCaseI}
Assume that 
\[
f:X_1\times \RR^m\rTo X_2\times\RR^m
\]
is as in \ref{SetupForThmIII}.
If $\partialmax X_1$ is irreducible and $\dim(X_1)=n\geq 2$ and if
some wall tree of $X_1$ is of type (I), then there is an isometry
$\bar f:X_1\rTo X_2$ with $(\bar f\times\id_{\RR^m})_*=f_*$.

\proof
If the wall tree $T_1$ is of type (I), then
$X_1=\EB (\partialmax X_1)$ by \ref{WallTreesAndProjectivities}.
Since $\partialmax X_1$ is thick and irreducible,
$\Pi_{\partialmax X_1}^+(a)$ acts $2$-transitively on
the ends of $T_1$ by \ref{KnarrThm}. 
By \ref{MainTreeProp}, the wall trees of
$X_2$ are also of type (I). Therefore
$X_2=\EB (\partialmax X_2)$. The isomorphism
$f_*:\partialmax X_1\rTo \partialmax X_1$ extends
to an isometry $\bar f:X_1\rTo X_2$ of the respective Euclidean cones.
\qed
\end{Lem}
\begin{Lem}
\label{IrredCaseIIandIII}
Assume that 
\[
f:X_1\times \RR^m\rTo X_2\times\RR^m
\]
is as in \ref{SetupForThmIII}.
If $\partialmax X_1$ is irreducible and $\dim(X_1)=n\geq 2$ and if
no wall tree of $X_1$ is of type (I), then the metric on $X_2$ can be
rescaled so that there is an isometry
\[
\bar f:X_1\rTo X_2
\]
with $(\bar f\times\id_{\RR^m})_*=f_*$.
This isometry $\bar f$ is unique and there is a constant
$r>0$ such that 
\[
d(f_1(x\times v),\bar f(x))\leq r\quad\text{ for all }\quad 
x\times v\in X_1\times\RR^m.
\]

\proof
Let $(a,b)$ be a pair of opposite panels in $\partialmax X_1$.
Then $(a,b)$ determines a \emph{wall} (a sphere of dimension $n-2$)
in the spherical building $\partialmax X_1$. 
$\partialmax X_1$ is irreducible, there are at most $2$
types of walls in $\partialmax X_1$.

The group $\Pi^+_{\partialmax X_1}(a)$
acts $2$-transitively on the ends of the wall tree $T_1$ by \ref{KnarrThm}.
We may apply \ref{MainTreeProp} to the coarse equivalence
$T_1\times\RR^{n-1+m}\rTo T_2\times\RR^{n-1+m}$, since
the equivariance condition is satisfied by our previous discussion.
Once and for all, we rescale the metric on $X_2$ in such
a way that $f_*$ extends to an equivariant isometry 
$\tau:T_1\rTo T_2$.

If $(a',b')$ is any other pair of panels 
in $\partialmax X_1$ of the same type as $(a,b)$, with
$X_1(a',b')=T_1'\times\RR^{n-1}$, 
then there is some projectivity which induces an isometry
$\phi_1:T_1\rTo T_1'$.
Pushing this projectivity forward via $f_*$, we obtain an
isometry $\phi_2:T_2\rTo T_2'$,
where $X_2(f_*a',f_*b')=T_2'\times\RR^{n-1}$.
By construction, the maps
$\phi_2\circ\tau\circ \phi_1^{-1}$
and $f_*$ induce the same map $\partial T_1'\rTo \partial T_2'$.
By \ref{MainTreeProp}, the map $\phi_2\circ\tau\circ \phi_1^{-1}$
is the unique equivariant tree isometry which accompanies
the coarse equivalence $X_1(a',b')\rTo X_2(f_*a',f_*b')$.
This shows that with respect to our metrics on $X_1$ and $X_2$,
the map $f_*:\partialmax X_1\rTo\partialmax X_2$
is ecological for all wall trees whose wall in
$\partialmax X_1$ is of the same type as the wall
determined by $(a,b)$.

Suppose that $z\in X_1$ is a thick point in an affine apartment
$A\subseteq X_1$. Since $X_1$ is
irreducible, $z$ is the intersection of $n$ walls 
$M_1,\ldots,M_n\subseteq A$ which are of the same type as the wall
determined by $(a,b)$. To see this, we note that the $W$-orbit
of any nonzero vector spans the ambient Euclidean space.
Each of these walls
determines a branch point in a panel tree. Let 
$M_{1,*},\ldots,M_{n,*}$ be the corresponding walls in $f_*A$,
defined by the isometries between the corresponding wall trees.
The intersection $M_{1,*}\cap\cdots\cap M_{n,*}$
is a point $z_*\in f_*A$. By \ref{MainTreeProp} there
is a uniform constant $s>0$ such that $d(f_1(z\times v),z_*)\leq s$.
If $z'\in A$ is another thick point and if $z'_*$ is constructed in the
same way, then the $n$ tree isometries yield $d(z,z')=d(z_*,z'_*)$,
whence $d(f_1(z\times v),f_1(z'\times v)\leq d(z,z')+2s$.
The thick points are cobounded in $X_1$, so the map
$x\mapstoo f_1(x\times v)$ is a rough isometry
$X_1\rTo X_2$. This implies by \ref{MainTreeProp}
that for no wall tree of $X_1$,
the accompanying isometry requires any rescaling of $X_2$.
The accompanying tree isometries therefore fit together
with $f_*$ to an ecological building isomorphism.
By Tits' result \ref{TitsEcologicalTheorem}, $f$ is accompanied by an isometry
$\bar f:X_1\rTo X_2$. This isometry maps the thick point
$z$ precisely to the point $z_*$ described above, $z_*=\bar f(z)$.
It follows that there is a constant $r>0$ such that
$d(f_1(x\times v),\bar f(x))\leq r$ holds for all
$x\times v\in X_1\times\RR^m$.
\qed
\end{Lem}
We now decompose the Euclidean building $X_1$ into a product
$X_1^I\times X_1^{II}\times X_1^{III}$ of Euclidean buildings of
types (I), (II) and (III), respectively. Similarly, we decompose
$X_2$. The next result implies Theorem III.
\begin{Thm}
\label{MainTheoremOfPaper}
Let 
\[
f:X_1\times\RR^{m_1}\rTo X_2\times\RR^{m_2}
\]
be a coarse equivalence of Euclidean buildings. Assume that
$\partialmax X_1$ and $\partialmax X_2$ are thick and that
$X_1$ splits off no tree factors. Then the following hold.

(o) There are numbers $m,n$ with $m_1=m_2=m$ and $\dim(X_1)=\dim(X_2)=n$.

(i)
The irreducible factors of $X_2$ can be rescaled in such a way that
there is an isometry 
\[
\bar f:X_1\rTo X_2\quad\text{ with }\quad f_*=(\bar f\times\id_{\RR^m})_*.
\]

(ii)
If 
\[
X_1=X_1^I\times X_1^{II}\times X_1^{III}\text{  and }
X_2=X_2^I\times X_2^{II}\times X_2^{III}
\]
are decomposed as above,
then $\bar f$ factors as a product,
\[
\bar f(x_I\times x_{II}\times x_{III})=
\bar f_I(x_I)\times\bar f_{II}(x_{II})\times\bar f_{III}(x_{III}).
\]
There is a constant $r>0$ such that
\[
d(\bar f_N(x_N),\pi_{X_N}f(x_I\times x_{II}\times x_{III}\times p))\leq r,
\]
for $N\in\{II,III\}$ and for all 
\[
x_I\times x_{II}\times x_{III}\times v
\in X_1^I\times X_1^{II}\times X_1^{III}\times\RR^m.
\]

(iii)
The maps $f_{II}$ and $f_{III}$ are unique.

\proof
We proceed by induction on the number irreducible factors of $X_1$.
The case of one irreducible factor is covered by \ref{IrredCaseI} and
\ref{IrredCaseIIandIII}.
In general, we decompose $X_1=Y_1\times Z_1$, with $Z_1$ irreducible.
As $f_*$ is an isomorphism, we have a corresponding decomposition
$X_2=Y_2\times Z_2$.
Fix opposite chambers $a,b$ in $\partialmax Y_1$, and let
$A\subseteq Y_1$ be the corresponding affine apartment.
Then $A\times Z_1$ is the union of all affine apartments
in $Y_1\times Z_1$ which contain $a,b$ at infinity. This relation
is preserved by $f_*$. So
if $y\times z\times v\in A\times Z_1\times\RR^m$, then
$f(y\times z\times v)$ has uniform distance from
$f_*A\times Z_1\times\RR^m$. It follows that 
$\pi_{f_*A\times Z_2\times\RR^m}\circ f$ is
a coarse equivalence between 
$A\times Z_1\times\RR^m$ and $f_*A\times Z_2\times\RR^m$.
By the induction hypothesis we can rescale the irreducible factors
of $Y_2$ in such a way that there is an isometry
between $Y_1$ and $Y_2$, and \ref{IrredCaseI} and
\ref{IrredCaseIIandIII} give us
isometries between $Z_1$ and $Z_2$, possibly after rescaling $Z_2$.
These isometries fit
together to an isometry $\bar f:X_1\rTo X_2$, with
$f_*=\bar f_*$. This gives (i). If $Z_1$ is of type (II) or
(III), then the claims (ii) and (iii) follow, by applying
\ref{IrredCaseIIandIII} to
$A\times Z_1\times\RR^m$ and $f_*A\times Z_2\times\RR^m$.
\qed
\end{Thm}


\renewcommand*{\thefootnote}{\fnsymbol{footnote}}
\setcounter{footnote}{0}
\setcounter{section}{1}
\setcounter{Thm}{0}
\renewcommand{\thesection}{\Alph{section}}

\newpage
\begin{center}
 \LARGE\bf Appendix to: `Coarse equivalences of Euclidean
buildings'\\
\end{center}
\begin{center}
by Jeroen Schillewaert\footnote{The first author is supported by Marie Curie IEF Grant GELATI (EC grant nr 328178)}
and Koen Struyve\footnote{The second author is supported
by 
 the Fund for Scientific Research --- Flanders (FWO - Vlaanderen)}
\end{center}

\bigskip\bigskip
The purpose of this appendix is to provide a generalization of the main results
of the preceding paper 
`Coarse equivalences of Euclidean buildings' by Linus Kramer and Richard Weiss.
We will refer to this paper as [KW]. 
The paper in question proves rigidity results for coarse equivalences of
metrically 
complete Euclidean buildings. We will show that the hypothesis of
completeness can be omitted.

The proof for this generalization is certainly not independent of the original
proof. In fact, we only 
discuss where and how the proof of Kramer and Weiss should be altered to allow
for non-complete Euclidean 
buildings. As our proof is an extension of the Kramer and Weiss argument, we
refer and rely on it for a 
detailed introduction, details and definitions. 

Here are the theorems that we obtain. 

\medskip\noindent\textbf{Theorem A.I\ }
{\em 
Let $G$ be a group acting isometrically on two leafless
trees $T_1,T_2$. Assume that there is a coarse equivalence
$f:T_1\rTo T_2$, that $T_1$ has at least $3$ ends
and that the induced map $\partial f:\partial T_1\rTo\partial T_2$
between the ends of the trees is $G$-equivariant. If the
$G$-action on $\partial T_1$ is $2$-transitive, then (after rescaling the
metric on $T_2$) there is a 
$G$-equivariant isometry $\bar f:T_1\rTo T_2$
with $\partial f=\partial \bar f$. If $T_1$ has at least two branch points,
then $\bar f$ is unique and has finite distance from $f$.}


\medskip\noindent\textbf{Theorem A.III\ }
{\em 
Let $f:X_1\times\RR^m\rTo X_2\times\RR^m$ be as in Theorem II in [KW] and
assume in addition that $X_1$ has no tree factors. 
Then there is (possibly after rescaling the metrics
on the de Rham factors of $X_2$) an isometry
$\bar f:X_1\rTo X_2$ with $(\bar f\times\id_{\RR^m})_*=f_*$.
Put $f(x\times v)=f_1(x\times v)\times f_2(x\times v)$.
If none of the de Rham factors of $X_1$ is a Euclidean cone over
its boundary, then $\bar f$ is unique and $d(f_1(x\times y),\bar f(x))$ is
bounded as a function of $x\in X_1$.}

\medskip\noindent
We now discuss the modifications one has to make in order to avoid the
assumption 
of metrical completeness in the results on coarse rigidity of Linus Kramer and
Richard Weiss in~[KW]. Completeness is used at two places in
their proof, which we will discuss separately.
All references are to~[KW] unless mentioned
otherwise.

\subsection*{Recovering the tree from the $G$-action}\label{section:1}
The first problem that occurs is that the Bruhat-Tits Fixed Point Theorem
requires
completeness, so a bounded isometry group acting on a metrically non-complete
$\RR$-tree $T$ does not necessarily have a fixed point (although it does in the
metric completion $\overline{T}$ of $T$). Consequently also 
Proposition~\ref{MaxBoundedSubgroups} no
longer holds.  The next lemma and propositions offer a way to deal with this
observation. For a geodesic segment $[x,y]$ in a tree we put $(x,y)=[x,y]\setminus\{x,y\}$.

\begin{Lem}\label{interior}
Let $T$ be an $\RR$-tree and $x,y \in \overline{T}$. Then the open segment
$(x,y)$ lies in $T$.

\proof
First of all note that the metric completion of an $\RR$-tree is again an 
$\RR$-tree by~\cite{Imr:77} (see also~\cite[Cor. II.1.10]{MorSha}), hence it
makes sense to speak about the open segment $(x,y)$. (Note however that the
metric
completion of a leafless metrically non-complete $\RR$-tree is never leafless.)

Given $\eps>0$, we choose points $x',y'\in T$ with $d(x,x'),d(y,y')<\eps$.
The geodesic segment $[x',y']$ is contained in $T$.
We put $x''=\pi_{[x',y']}(x)$ and $y''=\pi_{[x',y']}(y)$.
The unique geodesic from $x$ to $x'$ then passes through $x''$.
In particular,
$d(x,x'')<\eps$, and similarly $d(y,y'')<\eps$. By (T2), we have that
$[x,y]=[x,x'']\cup[x'',y'']\cup[y'',y]$.
The claim follows now, because
$[x'',y'']\subseteq T$ and because $\eps$ was arbitrarily small.
\qed
%
%
%
%
%
%
%
%
%
%
%
\end{Lem}
For the purpose of this appendix,
we restate condition (2-$\partial$) without the completeness assumption:

\medskip\noindent
{\boldmath\textbf{(2-$\partial$)}}
The group $G$ acts isometrically on the leafless tree $T$, and
this action is $2$-transitive on $\partial T$.

\medskip\noindent
The next two propositions serve as replacements of 
Propositions~\ref{AllPointsGIsolated} and 
\ref{MaxBoundedSubgroups} of~[KW]. 

\begin{Prop}
Assume  that the tree $T$ satisfies (2-$\partial$) and that the set of branch
points is dense. Then every point $x \in T$ is $G$-isolated in $\overline{T}$.

\proof
Let $x \in T$. Suppose that $G_x$ fixes another point $y \neq x$ in
$\overline{T}$. 
Let $y' \in (x,y)$. Then $y'$ lies in $T$ by Lemma~\ref{interior}  and $G_x$
fixes the geodesic segment $[x,y']$. 
There is a branch point $z$ between $x$ and $y'$, so
Proposition~\ref{TransThm} 
implies that
$y'$ is not a fixed point of $G_{x,z}$. \qed
\end{Prop}

This modification together with the Bruhat-Tits Fixed Point Theorem allows us,
as in~[KW], to conclude that the stabilizers of $G$-isolated
points in $T$ are maximal bounded subgroups.

\begin{Prop}
Assume (2-$\partial$) and that $P \subseteq G$ is a maximal bounded subgroup.
Then $P$ is the stabilizer of a $G$-isolated point, or it only fixes exactly one
point in the completion $\overline{T}$ and none in $T$.

\proof
Note that we may assume that the tree is metrically non-complete, and hence 
that the tree $T$ is of type (III) with a dense  set of branch points (see
Corollary~\ref{TreeStructureCor}).

Assume that what we want to prove is false.
As in Proposition~\ref{MaxBoundedSubgroups}, this implies that
$\overline{T}^P$ contains a geodesic segment $[x,y]\subseteq\overline{T}$ with $x \neq y$. Applying
Lemma~\ref{interior} one has that $(x,y) \subset T$, 
so there exists a subsegment $[x',y'] \subset (x,y) \subset T^P$ with $x' \neq
y'$. This situation was proved to be impossible in
Proposition~\ref{MaxBoundedSubgroups}.
\qed 
\end{Prop}
We are still left with the possibility that there are more maximal bounded
subgroups than only those corresponding to $G$-isolated points of $T$. However
those corresponding to $G$-isolated points can be recognized as follows:

\begin{Lem}
A maximal bounded subgroup $H$ is the stabilizer of a $G$-isolated point of $T$ 
if and only if either $T$ is not of type (III), or if $T$ is of type (III) and 
for each two ends $u,v \in \partial T$, 
$H$ is the only maximal bounded subgroup
of $G$ containing $\bra{(G_u \cup G_v) \cap H}$.

\proof
Let $H$ be a maximal bounded subgroup and let $x$ be the point stabilized by
$H$.
If the tree is not of type (III), then the tree is metrically complete and no
extra maximal bounded subgroups appear. 
If the tree is of type (III) then Lemma~\ref{RecoverApts} states that 
for each two ends $u,v \in \partial T$, the group
$H$ is the only maximal bounded subgroup
of $G$ containing $\bra{(G_u \cup G_v) \cap H }$ if and only if $x$ lies in some
apartment.
This holds if and only if $x$ is $G$-isolated.
\qed
\end{Lem}
Note that one can still distinguish between types (0)-(III) using
Lemma~\ref{interior} and \ref{Recover_d}.

Now we can modify the last part of Section 2 (starting directly after 
Proposition~\ref{MaxBoundedSubgroups}), and define $i_G(T)$ to correspond
with the set of maximal bounded subgroups of $G$ satisfying the conditions of
the above lemma. We conclude that Proposition~\ref{GroupEncodesTree} still
holds in the metrically non-complete case. 

Similarly one has to alter Proposition~\ref{BoundedThm} and
Theorem~\ref{MainTreeProp} restricting to the maximal bounded subgroups
with the above property. Observe that this property is preserved by the
$G$-equivariance.

\subsection*{\boldmath The use of retraction maps $\pi$}\label{section:2}
In \ref{CAT(k)Definition} one refers to~\cite[II.2.4]{BH} for a
$1$-Lipschitz retraction map $\pi_K: X \rTo K$ when $K$ is a metrically
complete convex subset in a CAT(0)-space $X$. This kind of retraction maps
is used in various places of~[KW]. 
In most of these places, $K$ is a closed convex subset of an apartment
and there is no problem. The only place where this is not the case is at the
beginning of \ref{SetupForThmIII}.

For each point $x$ in the completion of a set $K$ and for each $\eps > 0$, 
one can find a point $x' \in K$ such that $d(x,x') < \eps$. Combining this
with the above mentioned result for complete convex subsets of CAT(0)-spaces one
obtains the following lemma:

\begin{Lem}
Given a closed convex subset $K$ of a CAT(0)-space $X$ and an $\eps > 0$,
there exists a controlled map $\pi'_K : X \rTo K$ with control function
$\rho(t) = t +\eps$.
\end{Lem}
These maps $\pi'_K$ are sufficient for the purposes of 
\ref{SetupForThmIII}. Indeed, one only uses $1$-Lipschitz to obtain that
certain compositions of maps are controlled. This conclusion remains valid with
the weaker control function of $\pi'_K$.

\bigskip
{\raggedright
Linus Kramer\\
Mathematisches Institut, 
Universit\"at M\"unster,
Einsteinstr. 62,
48149 M\"unster,
Germany\\
e-mail: {\tt linus.kramer{@}uni-muenster.de}\\\smallskip
Jeroen Schillewaert\\
Department of Mathematics, Imperial College London, South Kensington Campus, London SW7~2AZ, England\\
e-mail: jschillewaert@gmail.com\\\smallskip
Koen Struyve\\
Department Mathematics,
Ghent University,
Krijgslaan 281, S22,
B-9000 Ghent,
Belgium\\
e-mail: {\tt kstruyve{@}cage.ugent.be}\\\smallskip
Richard M.~Weiss\\
Dept.~of Mathematics,
Tufts University,
503 Boston Ave.,
Medford, MA~02155, USA\\
e-mail: {\tt richard.weiss{@}tufts.edu}}

\end{document}